\documentclass[11pt]{amsart}
\usepackage{amsbsy,amsfonts}
\usepackage{latexsym,eucal,amsmath,amsthm}
\usepackage{amssymb}

\newcommand{\db}{\mathbb }

\newcommand{\scripte}{{\mathcal E}}
\newcommand{\scripto}{{\mathcal O}}

\newcommand{\Schwarz}{\mathcal S}

\def\Re{\mbox{\rm Re}}
\def\pc{\mbox{\rm pc}}

\def\H{{\bf H}}
\def\E{{\bf E}}
\def\Z{{\db Z}}

\def\R{{\db R}}
\def\C{{\db C}}
\def\torus{\mathbb T}
\def\integers{\mathbb Z}
\def\complex{\mathbb C}
\def\naturals{\mathbb N}

\def\rp{{^{-1}}}

\def\eps{\varepsilon}
\def\be#1{\begin{equation} \label{#1}}
\def\bas{\begin{align*}}

\newtheorem{theorem}{Theorem}
\theoremstyle{definition}

\theoremstyle{remark}
\newtheorem{remark}{Remark}
\theoremstyle{proposition}
\newtheorem{proposition}{Proposition}
\theoremstyle{lemma}
\newtheorem{lemma}{Lemma}
\theoremstyle{corollary}

\numberwithin{equation}{section}
\numberwithin{lemma}{section}
\numberwithin{remark}{section}

\begin{document}

\title[Ill-posedness for canonical defocusing equations]
{Asymptotics, frequency modulation, and low regularity ill-posedness for canonical defocusing equations}

\author{Michael Christ}
\thanks{M.C. is supported in part by N.S.F. grant DMS 9970660.}
\address{University of California, Berkeley}
\author{James Colliander}
\thanks{J.C. is supported in part by N.S.F. grant DMS 0100595 and N.S.E.R.C.
Grant RGPIN 250233-03.}
\address{University of Toronto}
\author{Terence Tao}
\thanks{T.T. is a Clay Prize Fellow and is supported in part by grants
from the Packard and Sloan Foundations.}
\address{University of California, Los Angeles}



\begin{abstract}
In a recent paper \cite{kpv:counter}, Kenig, Ponce and Vega study the low regularity behavior of the focusing nonlinear Schr\"odinger (NLS), focusing modified Korteweg-de Vries (mKdV), and complex Korteweg-de Vries (KdV) equations.  Using soliton and breather solutions, they demonstrate the lack of local well-posedness for these 
equations below their respective endpoint regularities.

In this paper, we study the defocusing analogues of these equations,  
namely defocusing NLS, defocusing mKdV, and real KdV, all in one spatial dimension, for which suitable soliton and breather solutions are unavailable. 
We construct for each of these equations 
classes of modified scattering solutions, 
which exist globally in time, and are asymptotic to solutions
of the corresponding linear equations up to explicit phase shifts.
These solutions are used to demonstrate lack of local well-posedness
in certain Sobolev spaces,
in the sense that the dependence of solutions upon initial data fails
to be uniformly continuous.
In particular, we show that the mKdV flow is not uniformly continuous in the 
$L^2$ topology, despite the existence of global weak solutions at this 
regularity.

Finally, we investigate the KdV equation at the endpoint regularity 
$H^{-3/4}$, and construct solutions for both the real and complex 
KdV equations. The construction provides a nontrivial time interval
$[-T,T]$
and a locally Lipschitz continuous map taking the initial data in
$H^{-3/4}$ to a distributional 
solution $u \in C^0 ([-T,T];
H^{-3/4})$ which is uniquely defined for all smooth data. 
The proof uses a generalized Miura transform 
to transfer the existing endpoint regularity theory for mKdV
to KdV. 
\end{abstract}

\date{1 July 2002}

\maketitle

\tableofcontents

\section{Introduction}

The purpose of this paper is to study asymptotic behaviour of solutions,
and low regularity well-posedness, for the defocusing 
nonlinear Schr\"odinger equation (NLS), the defocusing modified Korteweg-de Vries equation (mKdV), and the real Korteweg-de Vries equation (KdV).  This may be viewed as a follow-up to the work of Kenig, Ponce, and Vega \cite{kpv:counter} on the focusing
analogues of these equations.  We work on the real line $\R$; the case of the torus $\R/2\pi\Z$ is substantially easier and described in the last section.
In the next three sub-sections we describe our results for each of these equations 
in turn.

\subsection{The defocusing nonlinear Schr\"odinger equation}

The Cauchy problem for the cubic one-dimensional defocusing nonlinear Schr\"odinger equation\footnote{The minus sign in $-iu_t$ is convenient for our purposes, but can be removed if desired by replacing $u$ with $\overline{u}$.} (NLS) is
\be{nls}
\left\{
\begin{matrix}
-iu_t + u_{xx} = |u|^2 u; & u:[-T,T] \times \R_x \longmapsto \C,\\
u(0,x) = u_0(x),
\end{matrix}
\right.
\end{equation}
where $u_0$ is an element of a Sobolev space $H^s_x(\R)$ for some $s \in \R$.
The nonlinear Schr\"odinger equation is of widespread relevance in wave phenomena \cite{SulemBook}, \cite{NovikovText}. 
Indeed, the expectation is that whenever a physical system under consideration is described by a PDE which has a strongly dispersive linearization ($D^2_k \omega (k) \neq 0$), is weakly nonlinear and the solutions of interest are nearly monochromatic plane waves, NLS arises as an approximate model for the slowly varying wave amplitude. In this sense, NLS is a {\it{canonical}} dispersive
equation.

If the initial datum $u_0$ is in the Schwartz space $\Schwarz$, then there is a unique global smooth solution $u$ (see e.g. \cite{kato-kdv}).  In particular for each time $t$ we have a nonlinear 
evolution operator $S(t): \Schwarz \to \Schwarz$ defined by $S(t) u_0 := u(t)$, and a uniquely defined 
solution map $S: 
\Schwarz \to C^\infty_t(\R; \Schwarz)$ defined by $S u_0 := u$.

We are interested in the question of whether the solution map $S$ can be extended to rough initial data, such as data in the Sobolev space $H^s_x$ for some $s \in \R$.  If for every radius $R > 0$, there exists a time $T = T(R) > 0$ such that the solution map $S$ can be uniformly continuously and uniquely 
extended to a map from the ball $\{ u_0 \in H^s_x: \| u_0 \|_{H^s_x} <
R \}$ to the space\footnote{Of course, we endow the space
  $C^0_t((-T,T); H^s_x)$ with the topology induced by the norm
  $\sup_{t \in (-T,T)} \| u(t) \|_{H^s_x}$.} 
then 
we say that the equation \eqref{nls} is\footnote{
%
%
%
%
%
This ({\it{minimal}}) notion of local well-posedness is designed to
provide meaning to rough solutions obtained through a limiting
procedure of smooth functions. It differs subtly from what may be the
``most natural'' definition: For any $R>0$ there exists $T = T(R)>0$
such that the data-to-solution map $S$ is uniformly continuous and
uniquely defined from the ball $\{ u_0 \in H^s_x : {{\| u_0
    \|}_{H^s_x}} <R \}$ to the space $C^0 ( [-T, T] ; H^s_x )$. An
{\it{alternative}} notion of local well-posedness is also in common
use which replaces the space $C^0 ( [-T, T] ; H^s_x )$ in the ``most
natural'' definition by 
$C^0 ( [-T, T] ; H^s_x )\bigcap Y_T$ where $Y_T$ is an auxiliary
Banach space of functions of spacetime. This alternative
well-posedness is stronger than the ``most natural'' notion in the
sense that it provides extra $Y_T$-regularity of the solution but
weaker in the sense that the uniqueness property is in the smaller
intersected space. Of course, if the space $Y_T$ contains all smooth
solutions (which in practice it always does) then this alternative
local well-posedness implies the minimal notion of local
well-posedness defined in the main text. The positive results of
Kenig, Ponce and Vega we quote from
\cite{kpv:gkdv} and \cite{kpv:kdv} establish well-posedness with an
appropriate space $Y_T$ which contains all smooth solutions.}
{\it{locally well-posed}}.
%
%
If one can make $T$ arbitrarily large\footnote{This is not quite the same as setting $T=+\infty$, as the uniform continuity of the solution map may be destroyed in the infinite time limit.} and independent of $R$, then we say that \eqref{nls} is \emph{globally well-posed in $H^s_x$}.  
Since proofs of well-posedness based on a fixed point argument provide analytic dependence on the initial data, it is also natural to consider a more restricted notion of well-posedness requiring smoother dependence (e.g. $C^k$) upon the data than uniform continuity.

The following result is due to Tsutsumi:

\begin{theorem}\label{tsu}\cite{tsutsumi}  If $s \geq 0$, then the equation \eqref{nls} is globally well-posed in $H^s_x$.
\end{theorem}

This raises the question of what happens for $s < 0$.  The scale invariance
$$ u(t,x) \mapsto \frac{1}{\lambda} u\left ( \frac{t}{\lambda^2}, \frac{x}{\lambda} \right) $$
suggests that local well-posedness should fail for $s < -1/2$, while the Galilean invariance
$$ u(t,x) \mapsto e^{i\alpha x/2} e^{i\alpha ^2t/4} u(t, x + \alpha t)$$
suggests that local well-posedness should fail for $s < 0$.  This is because the spaces $\dot H^{-1/2}_x$ and $L^2_x$ are invariant under scaling, 
and under Galilean transformations, respectively.  However, these arguments are merely heuristic and do not constitute a rigorous proof of ill-posedness.

In \cite{kpv:counter}, Kenig, Ponce and Vega 
extended their earlier work \cite{BKSPV} with Birnir and Svanstedt (see also \cite{BPS}) and studied the \emph{focusing} analogue of \eqref{nls}, in which the nonlinear term $|u|^2 u$ is replaced by $-|u|^2 u$. All the above results for the defocusing equation extend to the focusing case, and furthermore there exist soliton solutions in this case.  By using the scale and Galilean invariances with these special soliton solutions, Kenig, Ponce and Vega show that the focusing NLS equation is not locally\footnote{Despite this negative result below $L^2$, it is still possible to obtain local well-posedness in certain spaces ``rougher'' than $L^2$ if one abandons the Sobolev scale of regularity.  See \cite{vargas-vega}.} well-posed in $H^s_x$ for any $s < 0$.
More precisely, they proved that the solution map $S$, restricted to initial
data in the Schwartz class,
fails to be uniformly continuous in the required norms.

Another instance where illposedness for a defocusing equation has
been established is a paper of Lebeau \cite{lebeau} on
the real-valued supercritical defocusing wave equation in $\R^3$ with
nonlinearity $u^7$. A more dramatic form of ill-posedness
is demonstrated there. Burq, G\'erard, and Tzvetkov \cite{BGT:NLSsphere}
have proved ill-posedness of the nonlinear Schr\"odinger equation \eqref{nls} in the
periodic case, and have obtained interesting ill-posedness results related
to global geometry on higher-dimensional spheres.

There are a number of papers in the literature in which it is shown
that the solution operator for various nonlinear equations fails
to be $C^k$ for some specific value of $k$, e.g. $k=1$ or $k=2$,
or fails to be Lipschitz continuous. 
Among the recent works along these lines are
\cite{brennerkumlin}, \cite{ntt}, \cite{tzvetkov}.
To establish such results
amounts to showing that certain multilinear operators acting on Sobolev
spaces fail to satisfy certain {\it a priori} inequalities.
In contrast, ill-posedness results like those
of \cite{kpv:counter} and \cite{lebeau} depend on the
analysis of certain exact {\em solutions} of the equations.  
The fundamental issue in the defocusing case
is how to construct a suitable family of solutions,
which are sufficiently ill-behaved to demonstrate ill-posedness.

The first result of this paper extends the results of \cite{kpv:counter}
concerning NLS to the defocusing case.

\begin{theorem}\label{defocus-nls}  The Cauchy problem \eqref{nls} is not locally well-posed in $H^s_x$ for any $s < 0$.
\end{theorem}

We prove this theorem in Sections \ref{pseudo-sec}, \ref{nls-sec}, again by showing
that the restriction of $S$ to Schwartz class initial data fails to be uniformly
continuous.  (This map was previously shown in \cite{borg:xsb} in both focusing and defocusing cases to be non-analytic, and indeed not even $C^2$ in $H^s_x$ for $s < 0$.  Similar results were also obtained for the mKdV and KdV equation in \cite{borg:xsb}, \cite{tzvetkov}.)  The new difficulty in this case is the lack of soliton solutions (or indeed of any non-zero exact solution).  An inspection of the arguments in \cite{kpv:counter} reveals that soliton/breather solutions are not essential; instead, it suffices to have
smooth solutions to \eqref{nls} whose global time development can be accurately controlled,
and is sufficiently sensitive to initial conditions.
The scale and Galilean invariances can then be used to convert such solutions to a family of solutions whose dependence on their initial data can be made arbitrarily
non-uniform in $H^s_x$ norm for any $s < 0$. Observe that this scaling procedure 
can convert long times to arbitrarily short times, which explains our desire for 
global control of solutions.

To construct such global smooth solutions we shall use the modified scattering
asymptotics 
introduced by Ozawa\footnote{The authors are indebted to Kenji Nakanishi for suggesting this approach.} \cite{ozawa}.  Following Ozawa, we first use the pseudo-conformal transformation to convert the global problem to a local one, and then approximate the PDE by 
an associated ODE, using energy methods to estimate the error.  This approach will be expanded upon in detail in Section \ref{pseudo-sec}, with the application to ill-posedness given in Section \ref{nls-sec}. 
It works equally well for focusing NLS.  

Interestingly, our arguments are different in the subcritical $s>-1/2$, critical $s=1/2$, and supercritical $s<-1/2$ cases.
For the supercritical analysis, we introduce a different
construction of solutions, based on an approximation by the
zero dispersion limit $-iu_t = |u|^2 u$ of the equation. 
These solutions are controlled only for short times, yet still suffice
for the ill-posedness argument.
After this paper was nearly completed, we discovered how a
frequency-modulated version of this zero dispersion limit
construction could be used to obtain ill-posedness of NLS
for arbitrary $s<0$; thus scattering-type solutions could be
eliminated from the discussion. We plan to discuss this in
a subsequent paper. 

\subsection{The defocusing modified Korteweg-de Vries equation}

The Cauchy problem for the defocusing modified Korteweg-de Vries equation is
\begin{equation}\label{mkdv}
\left\{
\begin{matrix}
u_t + u_{xxx} = 6 u^2 u_x\\
u(0,x) = u_0(x),
\end{matrix}
\right.
\end{equation}
where $u$ is a real\footnote{The distinction between focusing and defocusing caes
is only meaningful in the real setting here, because substituting $u=iv$
converts one to the other. 
Note also that the global well-posedness result in Theorem 3 is only known to be valid in the real-valued setting.}-valued function on $\R \times (-T,T)$ and $u_0 \in H^s_x(\R)$ for some $s \in \R$.  
The modified KdV equation arises as a natural extension of the 1d wave equation \cite{NovikovText} and is therefore a canonical dispersive equation. As with NLS, one has unique global smooth solutions from Schwartz data 
(see e.g.\ \cite{kato}), and so again we have a solution operator $S$.  We can then ask whether \eqref{mkdv} is locally well-posed in $H^s_x$.

This question was studied in a series of papers \cite{kato-kdv}, \cite{SautThesis}, \cite{SautJDE79}, \cite{KatoManuscripta}, \cite{kato}, \cite{gtv:gkdv}, \cite{kpv:gkdv}, \cite{ckstt:2} culminating in the following local and global results.

\begin{theorem}\label{mkdv-lwp} The equation \eqref{mkdv} is locally well-posed in $H^s_x$ for $s \geq 1/4$ \cite{kpv:gkdv}, and is globally well-posed in $H^s_x$ for $s > 1/4$ \cite{ckstt:2}.
\end{theorem}

It is likely that one also has global well-posedness at the endpoint $s=1/4$, 
but we do not pursue this question here.  In Section \ref{mkdv-lwp-sec} we review the arguments from \cite{kpv:gkdv} as we shall need them again here\footnote{It is a paradoxical fact that, in the absence of exact solutions, we need the well-posedness theory at high regularities in order to prove the ill-posedness at low regularities.  This is because we need some sort of well-posedness to control the low regularity solution accurately enough to quantify the ill-posedness.  See also Lemma \ref{modulate}, in which smooth functions are used to construct rough ones.}.

As with NLS, there is a scale invariance
\be{mkdv-scaling}
u(t,x) \mapsto \frac{1}{\lambda} u(\frac{t}{\lambda^3}, \frac{x}{\lambda})
\end{equation}
which again suggests ill-posedness for $s < -1/2$, however there is no exact analogue of the Galilean invariance of NLS.  On the other hand, we have (formally at least) a conservation law for the $L^2$ norm
\be{l2-cons}
\int u(T,x)^2\ dx = \int u_0(x)^2\ dx \hbox{ for all } T \in \R
\end{equation}
so one might hope to have some sort of local well-posedness at the $L^2$ level.

In complete analogy with NLS, there is a focusing variant of \eqref{mkdv} in which the nonlinearity $6 u^2 u_x$ is replaced by $-6 u^2 u_x$, and Theorem \ref{mkdv-lwp} extends to the focusing case.  As with NLS, the focusing case admits a rich family of soliton
solutions, and in \cite{kpv:counter} Kenig, Ponce and Vega were able to show failure of local well-posedness in $H^s$ for $s < 1/4$, despite the conservation
law \eqref{l2-cons}.   
No solitons are available for the defocusing mKdV equation.
Nevertheless, we extend this ill-posedness result to the defocusing case:

\begin{theorem}\label{defocus-mkdv} The Cauchy problem \eqref{mkdv} is not locally well-posed in $H^s_x$ for any $-1/4 < s < 1/4$; more precisely, the solution operator fails
to be uniformly continuous with respect to the $H^s$ norm.
\end{theorem}

It seems likely that the restriction $s>-1/4$ is merely an artifact
of our method. Perhaps arguments related to the alternative we use
in the supercritical range $s<-1/2$ for the nonlinear Schr\"odinger equation
might be extended and adapted to this case.

It is interesting to view this theorem at $s=0$ in the light of
Kato's construction \cite{kato} of global weak solutions in $L^2_x$.  
Thus the mKdV flow can be defined in $L^2_x$ in a weak sense, but the resulting
flow is not uniformly continuous. It is natural to ask whether it
might be non-unique. 
These issues may be related to the nonuniqueness of weak solutions of the Euler equation \cite{SchefferNonunique}, \cite{Schnirelman}.

\eqref{mkdv} is defocusing in the sense that the positive definite quantity
$\int \big(u_x(x,t)^2+u(x,t)^4\big)\,dx$
is formally conserved; for focusing 
mKdV the corresponding conserved
quantity is not semidefinite.
Moreover, \eqref{mkdv} is also defocusing in a second sense: in certain 
asymptotic regimes, it can be approximated well by defocusing NLS as described
in the following paragraphs. 
See \cite{schneider}, \cite{BoydChen} for a related approximation of real KdV by NLS.

We prove Theorem \ref{defocus-mkdv} in Section \ref{crude-sec}, after a review of the local well-posedness theory in Section \ref{mkdv-lwp-sec}.  One idea is to approximate the mKdV equation by the NLS equation. 
In a different asymptotic regime, such an approximation has been carried out by Schneider \cite{schneider}.  We also give an alternate argument in Sections \ref{energy-sec} and \ref{approx-sec}.

Here we briefly sketch the means of approximation.  
Define the spacetime Fourier transform
$$ \tilde u(\tau,\xi) := \int e^{-it\tau} e^{-ix\xi} u(t,x)\ dt dx.$$
We expect time-localized 
solutions to the NLS equation to have spacetime Fourier transform near the parabola $\tau = \xi^2$, while time-localized solutions to the mKdV equation should have spacetime Fourier transform near the cubic $\tau = \xi^3$.

Pick a large number $N \gg 1$.  If we make the linear change of variables
$$ \tau = N^3 + 3^{1/2} N^{3/2} \xi' + \tau'; \quad \xi = N + 3^{-1/2} N^{-1/2} \xi'$$
then the cubic $\tau = \xi^3$ becomes
$$ \tau' = (\xi')^2 + 3^{-3/2} N^{-3/2} (\xi')^3.$$
Thus, for $|\xi'| \ll N$, this linear transformation converts the cubic to an approximate parabola.

Unraveling this transformation using the spacetime Fourier transform, we are led to the following heuristic: if $u$ solves the NLS equation \eqref{nls}, then the function
\be{approx}
v(t,x) := \sqrt{\frac{2}{3N}} \Re~  e^{iNx} e^{iN^3 t} u\left( t,\frac{x + 3N^2t}{3^{1/2} N^{1/2}} \right)
\end{equation}
is an approximate solution to \eqref{mkdv}.  
(To approximately  equate the cubic NLS nonlinearity $|u|^2 u$ with the cubic mKdV nonlinearity $6 u^2 u_x$ requires further calculations  which are 
omitted here;
those calculations give rise to the
factor $\sqrt{\frac{2}{3N}}$ and the real part operator $\Re$.
 The derivative in the mKdV nonlinearity is approximated
by $iN$ on the Fourier side, so that no derivative appears in the NLS approximation).  
If the NLS solution $u$ is mostly supported in the frequency range $|\xi'| \ll N$, then we observe that the $H^{1/4}$ norm of $v(t)$ is comparable to the $L^2$ norm of $u(t)$.
Thus Theorem \ref{defocus-mkdv} is closely related to Theorem \ref{defocus-nls}.

We are left with the problem of proving the existence of exact solutions of defocusing
mKdV which are well approximated by these solutions of NLS. This part of the analysis
is somewhat technical; we present two separate methods (one based on 
local smoothing and Strichartz estimates, one on energy estimates) 
for doing so and controlling the error in the approximation.  The method based on energy methods seems quite general, and should be able to yield a class of global solutions to a variety of equations.

\subsection{The Korteweg-de Vries equation}

Our final results concern the Cauchy problem for the Korteweg-de Vries equation (KdV)
\begin{equation}\label{kdv}
\left\{
\begin{matrix}
u_t + u_{xxx} = 6 u u_x \\
u(t,x) = u_0(x),
\end{matrix}
\right.
\end{equation}
where $u(t,x)$ is defined on $(-T,T) \times \R$ and is either real or complex-valued. 
The real-valued KdV equation arises as an approximate model to the standard real-valued 1d wave equation providing corrections due to weak nonlinearity and dispersion \cite{SulemBook}, \cite{NovikovText}. The universal relevance of the wave equation and these corrections justifies referring to the KdV as a canonical dispersive equation.    
Again (see e.g.\ \cite{kato-kdv}) smooth solutions of \eqref{kdv} are known to exist for Schwartz initial data.  Indeed, we have

\begin{theorem}\label{kdv-lwp}  In both the real and complex cases, the Cauchy problem \eqref{kdv} is locally well-posed in $H^{s}$ for $s > -3/4$ \cite{kpv:kdv}, and, in the real valued case, is also globally well-posed in $H^s$ for $s > -3/4$ \cite{ckstt:2}.
\end{theorem}

Our first result (proven in Section \ref{general-sec}) is to extend part of the local result to the endpoint $s=-3/4$.

\begin{theorem}\label{kdv-endpoint}  
In both the real and complex cases, for 
any $\phi \in H^{-3/4}$, there exist  $~T = T( \| \phi
\|_{H^{-3/4}} ) >0$ and a locally Lipschitz map
\begin{equation*}
  H^{-3/4} \ni \phi \longmapsto u \in C([0,T];
  H^{-3/4} )
\end{equation*}
which extends the smooth-data-to-solution map for
the initial value problem \eqref{kdv}. Moreover, for each $\phi \in
H^{-3/4}, ~u(t,x)$ is a weak solution of the KdV equation
\eqref{kdv}. 
\end{theorem}

Thus the initial value problem \eqref{kdv} is well-posed in
$H^{-3/4}$, in the minimal sense defined above. 

It is likely that one also has a global-in-time result at
$H^{-3/4}$ for the real KdV equation, but again we do not address these issues here.

A different endpoint result (with $H^{-3/4}$ replaced by a Besov variant) has been independently obtained recently by Muramatu and Taoka \cite{MuramatuTaoka}.  The bilinear estimate used to obtain the local results for $s > -3/4$ fails at the endpoint $s=-3/4$ \cite{ntt}; instead, we study
the endpoint $s = -3/4$ from the theory of the mKdV equation at $s=1/4$ by using a variant of the \emph{Miura transform} $u \mapsto u_x + u^2$, which maps solutions of defocusing mKdV to real KdV.\footnote{
The variant $u \mapsto u_x + iu^2$ maps focusing mKdV to complex KdV.}  Observe that this transform maps $H^{1/4}$ continuously to $H^{-3/4}$.  Unfortunately the Miura transform is not invertible; however, we will modify the Miura transform slightly (in a manner reminiscent of Gardner's 
modification of the Miura transform \cite{MiuraTransform}, see also \cite{MiuraKdVSurvey}) to make the transform invertible and close the argument. 

Next, we address the situation when $s < -3/4$.  The scale invariance
$$ u(t,x) \mapsto \frac{1}{\lambda^2} u\left( \frac{t}{\lambda}, \frac{x}{\lambda}\right)$$
suggests there is ill-posedness for $s < -3/2$.  Again there is no direct analogue{\footnote{The Galilean invariance $u(x,t) \longmapsto u(x + 6 \beta t, t) + \beta$ of KdV does not preseve decay properties at spatial infinity.}}
 of Galilean invariance, nevertheless there are breather solutions for complex KdV, and in \cite{kpv:counter} it was shown that the complex KdV equation is not locally well-posed in $H^s$ for any $s < -3/4$. 

In analogy with our prior results, we extend this result to the real case.

\begin{theorem}\label{kdv-illp}  The real KdV equation is not locally well-posed in $H^s$ for any $-1 \leq s < -3/4$; more precisely, the solution operator fails to
be uniformly continuous with respect to the $H^s$ norm.
\end{theorem}

This will be proved as a simple consequence of Theorem \ref{defocus-mkdv} and the Miura transform in Section \ref{miura-sec}.  The condition $-1 \leq s$ can certainly be relaxed but we do not pursue this matter here.

\subsection{Periodic analogues}

These results have analogues in the periodic case. Let
$\torus = \R/2\pi\integers$, and consider the same
partial differential equations now for $(t,x)\in\R\times\torus$.

\begin{theorem} \label{periodic-thm}
The defocusing nonlinear Schr\"odinger
equation is illposed in $H^s(\torus)$ for all $s<0$.

The modified real Korteweg-de Vries equation
is illposed in $H^s(\torus)$ for all $s\in(-1,1/2)$.

The real Korteweg-de Vries equation
is illposed in $H^s(\torus)$ for all $s\in(-2,-1/2)$.
\end{theorem}

Local and global well-posedness are known to hold for all larger
exponents $s$ \cite{borg:xsb}, \cite{kpv:kdv}, \cite{ckstt:gKdV},
\cite{ckstt:2}. As we shall see, the ill-posedness results are
substantially easier to obtain in the periodic case, although of the
same basic flavor. The first of these three conclusions has already  
been obtained by Burq, G\'erard and Tzvetkov \cite{BGT:NLSsphere}.

We believe that the
lower bounds on $s$ in the theorems for the KdV and mKdV equations
are merely artifacts of the method of proof.

\section{Notation and modulation bounds}\label{notation-sec}

$C$ denotes various constants depending only on $s$.  The notations
$A \lesssim B$ or $A = O(B)$ denote the estimate $A \leq CB$.

We define the spatial Fourier transform by
$$ \hat u(t,\xi) := \int e^{-ix \xi} u(t,x)\ dx.$$
The operator $\partial_x$ is conjugated to the multiplier $i\xi$ by the Fourier transform.

The following lemma will be used to estimate $H^s$ norms of
high-frequency modulations of smooth functions. 

\begin{lemma}\label{modulate}  Let $-1/2 < s$, $\sigma\in\R^+$
 and $u\in H^{\sigma}(\R)$. 
For any 
$M\ge 1$, $\tau\in\R^+ $, $x_0 \in \R$, and $A > 0$ let
$$ v(x)=v_{M,\tau,x_0,A}(x) := A e^{iMx} u((x-x_0)/\tau).$$

\noindent
(i) Suppose $s\ge 0$. 
Then there exists a constant $C_1<\infty$, depending only on $s$, 
such that whenever $M\cdot\tau\ge 1$,
$$ \| v \|_{H^s} \le C_1|A| \tau^{1/2} M^s \|u\|_{H^{s}}$$
for all $u,A,x_0$. 

\noindent
(ii) Suppose that $s<0$ and that $\sigma\ge |s|$.
Then there exists a constant $C_1<\infty$, depending only on $s$ and on $\sigma$, 
such that whenever $1\le \tau\cdot M^{1+(s/\sigma)}$,
$$ \| v \|_{H^s} \le C_1|A| \tau^{1/2} M^s \|u\|_{H^{\sigma}}$$
for all $u,A,x_0$.

%
%

\noindent
(iii) There exists $c_1>0$ such that
for each $u$ there exists $C_u<\infty$ such that
$$ \| v \|_{H^s} \ge c_1 |A| \tau^{1/2} M^s\|u\|_{L^2}$$
whenever $\tau\cdot M\ge C_u$.

\end{lemma}

\begin{proof}
\begin{align*}
A^{-2}\tau^{-1}M^{-2s}\|v\|_{H^s}^2
&= c\tau^{-1}M^{-2s}
\int (1+|\xi|^2)^{s} \tau^2 |\hat u(\tau(\xi-M))|^2\,d\xi
\\
&= 
c\int \left(\frac{\tau^2 + |M\tau+\eta|^2}{\tau^2 M^2}\right)^{s} 
|\hat u(\eta)|^2\ d\eta
\\
&\lesssim \int_{|\eta|\le \tau M/2}  |\hat u(\eta)|^2
+ 
\int_{\frac12\tau M\le|\eta|\le 2\tau M} M^{-2s} |\hat u(\eta)|^2
+ 
\int_{|\eta|\ge \tau M/2} \frac{|\eta|^{2s}}{(\tau M)^{2s}}|\hat u(\eta)|^2
\\& = I + II + III. 
\end{align*}
Term $I$ is $\lesssim \|u\|_{L^2}^2$. 
If $s\ge 0$ then $M^{-2s}\le 1$, so $II\lesssim \|u\|_{L^2}^2$, 
and $III\lesssim \|u\|_{H^s}^2$, because $\tau M\ge 1$.

If $s<0$ then $III \lesssim \|u\|_{L^2}^2$, since 
$|\eta|/\tau M\gtrsim 1$.
Moreover,
$II \lesssim M^{-2s}(\tau M)^{-2\sigma} \|u\|_{H^\sigma}^2$, 
which is $\lesssim \|u\|_{H^\sigma}^2$ under the further hypothesis
$1\le \tau\cdot M^{1+(s/\sigma)}$.


To obtain (iii), it suffices to consider term $I$:
for any $u$, $\int_{|\eta|\le \tau M/2} |\hat u(\eta)|^2$ approaches
$c\|u\|_{L^2}^2$ as $\tau M\to\infty$.
\end{proof}

\section{The pseudo-conformal transformation}\label{pseudo-sec}


\subsection{Definition of the $pc$ transform and some basic properties}

In this section we introduce the pseudo-conformal change of variables
\begin{equation}
\label{pseudoconformalchange}
(y,s) := \left( \frac{x}{t+1}, \frac{1}{t+1}\right); \quad (x,t) = \left(\frac{y}{s}, \frac{1}{s}-1\right)
\end{equation}
to analyze the asymptotic behavior of the NLS equation \eqref{nls} as $t \to +\infty$.  The $+1$ shift in time is 
introduced purely to avoid an artificial singularity at the initial time $t=0$ and should be ignored on a first reading.

A standard stationary phase computation (see e.g.\ \cite{heron}; alternatively, one can use the fundamental solution $C t^{-1/2} e^{-ix^2/4t}$) shows that solutions to the free Schr\"odinger equation $-iu_t + u_{xx} = 0$ behave asymptotically as $t \to +\infty$ like
$$ u(t,x) \approx (1+t)^{-1/2} \exp(-ix^2/4(t+1)) \phi(y)$$
for some function $\phi$ (which is essentially the Fourier transform of $u(-1)$).  

Motivated by this, we introduce the pseudo-conformal transformation $v = \pc(u)$, $u = \pc^{-1}(v)$ defined by the formulae
\begin{align}
u(t,x) &:= (1+t)^{-1/2} \exp(-ix^2/4(t+1)) v(s,y) \label{pc} \\
v(s,y) &:= s^{-1/2} \exp(iy^2/4s) u(t,x) \label{pc-inverse},
\end{align}
where it is understand that $(y,s)$ is always related to $(x,t)$ by 
the pseudo-conformal transformation \eqref{pseudoconformalchange}.

For each fixed time $t$, the map $u(t) \to v(s)$ is a linear isometry on $L^2$.  This map is not so well behaved on other Sobolev spaces $H^k$ because of the
highly oscillatory factors $e^{-ix^2/4(t+1)}$, $e^{iy^2/4s}$.  
To get around this we shall work in weighted Sobolev spaces, which we now discuss. 

For any integer $k \geq 0$, we define the space $\H^{k,k}_x$ to be the closure of Schwartz functions under the norm
$$ \| u \|_{\H^{k,k}_x} := \sum_{i,j \geq 0:\ i+j \leq k} \| x^i \partial_x^j u \|_{L^2_x}.$$
Thus the $\H^{k,k}_x$ norm controls the $H^k_x$ norm but also incorporates some spatial decay.  Roughly speaking, $\H^{k,k}_x$ is to the Hermite operator $-\Delta + |x|^2$ as $H^k_x$ is to the Laplacian $-\Delta$.

The next three simple lemmas control the behaviour of the $\H^{k,k}_x$ spaces under the pseudo-conformal transformation, pointwise multiplication, and the free Schr\"odinger flow.  To simplify the notation we shall often omit the variable $x$ from the norms $\H^{k,k}_x$, $H^k_x$, etc. when it is clear from context what the variable is.

\begin{lemma}\label{pc-convert}  Let $v = \pc(u)$.  Then 
$$ \| u(t) \|_{H^{k}_y} \lesssim \| v(s) \|_{\H^{k,k}_y}$$
for all $t \geq 0$ and integer $k \geq 0$, where the implicit constant depends on $k$ but not on $t$.  In the case $t=0$ we can improve this to
$$ \| u(0) \|_{\H^{k,k}_x} \lesssim \| v(1) \|_{\H^{k,k}_y}.$$
\end{lemma}

\begin{proof}
A brute force induction shows that a derivative $\partial_x^a u(t)$ of $u(t)$ can be expressed as a finite linear combination of terms of the form
$$ (1+t)^{-b-1/2} (\frac{x}{1+t})^c e^{-ix^2/4(t+1)} \partial_y^d v(s,y)$$
where $b, c, d$ are non-negative integers such that $c+d \leq a$.  The first claim follows.  The second claim follows from the identity $u(0,x) = \exp(-ix^2/4) v(0,x)$ and another brute force induction.
\end{proof}

\begin{lemma}\label{algebra}  For any $k \geq 1$, 
\be{a-1}
\| u v \|_{\H^{k,k}} \lesssim \| u \|_{\H^{k,k}} \| v \|_{\H^{k,k}}
\end{equation}
where the implicit constant is allowed to depend on $k$.  More precisely,
\be{a-2}
\| u v \|_{\H^{k,k}} \lesssim \| u \|_{\H^{k,k}} \| v \|_{L^\infty}
+ \| u \|_{\H^{k-1,k-1}} \| v \|_{\H^{k+1,k+1}}.
\end{equation}
If $u$ is real, then 
\be{a-3}
\| \exp(iu) v\|_{\H^{k,k}} \lesssim (1 + \| u \|_{\H^{k,k}})^k \| v \|_{\H^{k,k}}.
\end{equation}

\end{lemma}

\begin{proof}  
Let $i, j \geq 0$ be such that $i+j \leq k$.  Observe that $x^i \partial_x^j (uv)$ can be written as a finite linear combination of terms of the form $x^i (\partial_x^l u) (\partial_x^m v)$, where $l+m = j$.  At least one of $l,m$ must be less than or equal to $k-1$; without loss of generality we may assume $m \leq k-1$.  But then by Sobolev embedding we have $\| \partial_x^m v \|_\infty 
\lesssim \| v \|_{H^k} \lesssim \| v \|_{\H^{k,k}}$.  The claim \eqref{a-1} then follows by H\"older's inequality.

To prove \eqref{a-2} we refine the above argument.  If $l=k$ then we can take $v$ out in $L^\infty$ to estimate this term by $\| u \|_{\H^{k,k}} \| v \|_{L^\infty}$.  If $l < k$ then we take $\partial_x^m v$ out in $L^\infty$ and use Sobolev to majorize
this term by $\| u \|_{\H^{k-1,k-1}} \| v \|_{\H^{k+1,k+1}}$.

The inequality \eqref{a-3} is proven similarly to \eqref{a-1}, but one uses the chain rule $l$ times to expand out $\partial_x^l \exp(iu)$, and then discards the bounded factor $\exp(iu)$.  The details are left to the reader.
\end{proof}

\begin{lemma}\label{evolve}  If $t = O(1)$ and 
$1 \leq k \in {\mathbb{N}}$ then we have the estimate
$$ \| \exp(it \partial_{xx}) u \|_{\H^{k,k}_x} \lesssim \| u \|_{\H^{k,k}_x}$$
where the implicit constant depends on $k$.
\end{lemma}

In phase space, this Lemma asserts that the norm $|x| + |\xi|$ is stable under the flow $(x,\xi) \mapsto (x + t\xi, \xi)$ for $t = O(1)$.

\begin{proof}
Taking Fourier transforms, and observing from Plancherel and the product rule that $\| \hat u \|_{\H^{k,k}_x} \lesssim \| u \|_{\H^{k,k}_x}$, we see it suffices to show that
$$ \| \partial_\xi^l (\xi^m e^{i t \xi^2} \hat u) \|_{L^2_\xi} \lesssim \| \hat u \|_{\H^{k,k}_x}$$
whenever $0 \leq l,m$ and $l+m \leq k$.  By the product and chain rule, we can expand the left-hand side as a bounded linear combination (for $t = O(1)$) of terms of the form $\xi^a e^{i t\xi^2} \partial_\xi^b \hat u$, where $0 \leq a,b$ and $a+b \leq k$.  The claim follows.
\end{proof}

\subsection{Relation to Schr\"odinger equations}

We now return to nonlinear Schr\"odinger equations.
Some tedious computation using \eqref{pc-inverse} yields
\bas
v_s &:= s^{-5/2} e^{iy^2/4s} (-\frac{s}{2}u - \frac{iy^2}{4}u - y u_x - u_t) \\
v_y &:= s^{-5/2} e^{iy^2/4s} (\frac{iys}{2}u + s u_x) \\
v_{yy} &:= s^{-5/2} e^{iy^2/4s} (-\frac{y^2}{4}u + iyu_x + \frac{is}{2} u + u_{xx}) \\
\end{align*}
so that we have the identity
\be{pc-ident}
iv_s + v_{yy} = s^{-5/2} \exp(iy^2/4s) (-iu_t + u_{xx})
\end{equation}
for arbitrary $u$.  Because of this, the map $\pc$ transforms the Cauchy problem \eqref{nls} to a backwards Cauchy problem
\begin{equation}\label{back}
\left\{
\begin{matrix}
iv_s + v_{yy} = s^{-1} |v|^2 v \\
v(1,y) = v_1(y), & 0 < s \leq 1
\end{matrix}
\right.
\end{equation}
where $v_1(y) := e^{iy^2} u_0(y)$.

The singular term $1/s$ suggests that solutions of this equation should become singular in some sense as $s \to 0^+$.  Indeed, dropping the dispersive term $v_{yy}$ leaves the associated ODE\footnote{The authors are indebted to Kenji Nakanishi for the idea of introducing this ODE as an approximating equation.}  
\begin{equation}\label{ODE}
iv_s = s^{-1} |v|^2 v,
\end{equation}
for which there are explicit solutions $v = v^{[w]}$ of the form
\begin{equation}\label{explicit}
v^{[w]}(s,y) := w(y) \exp(-i |w(y)|^2 \log s)
\end{equation}
for any function $w(y)$. Now $v^{[w]}$ is singular, in the sense that it has
no limit as $s \to 0^+$.  
Moreover, for smooth $w$, $\partial_x^2 v^{[w]} \thicksim ( \log s )^2$
with the implicit constant depending upon $w$. Since $(\log s )^2 \ll s^{-1},$ we appear to be justified in ignoring the dispersive term in \eqref{back}.
This idea seems to be due to Ozawa \cite{ozawa}. 

We now show the following asymptotic completeness result, which is crucial to all of our ill-posedness results. 

\begin{lemma}\label{sing}  Let $K \geq 5$ be an integer, and let $w \in \H^{K+2,K+2}(\R)$ have an $\H^{K+2,K+2}$ norm of $O(\eps)$ for some small constant $0 < \eps \ll 1$.  Then if $\eps$ is sufficiently small, 
there exists $v_1 \in \H^{K,K}(\R)$ such that the unique solution
$v = v^{\langle w \rangle}$ to the backwards Cauchy problem \eqref{back} 
with initial datum $v_1$ satisfies
\be{asym}
\| v^{\langle w \rangle}(s) - v^{[w]}(s) \|_{\H^{K,K}} \lesssim \eps s (1 + |\log s|)^C \hbox{ for all } 0 < s \leq 1
\end{equation}
Furthermore, the map $w \mapsto v$ is Lipschitz continuous from the ball $\{ w \in \H^{K+2,K+2}: \| w \|_{\H^{K+2,K+2}} \leq \eps \}$ to $L^\infty_t((0,1]; \H^{K,K})$, i.e.
\be{lip}
\sup_{0 < s \leq 1} \| v^{\langle w' \rangle}(s) - v^{\langle w \rangle}(s) \|_{\H^{K,K}} \lesssim \| w' - w \|_{\H^{K+2,K+2}}
\end{equation}
for all $w$, $w'$ in the above ball.

\end{lemma}

\noindent
{\em Notational convention.\/}
Throughout the paper, square bracket superscripts (as in $v^{[w]}$) are used
to denote explicit, but approximate, solutions to nonlinear PDE, 
whereas angular bracket superscripts (as in $v^{\langle w \rangle}$) denote 
related exact solutions.  The functions $w$ will be in some sense ``data'' for these solutions, though not always in the classical sense of initial data.  Our strategy throughout is to first find an approximate solution $v^{[w]}$, then to carry out a perturbation analysis to pass from $v^{[w]}$ to an exact solution $v^{\langle w \rangle}$.

One can relax the condition $K \geq 5$ substantially, but this has no 
advantage here, and in fact for our ill-posedness application it will be useful to have $K$ arbitrarily large. 

\begin{proof}[Proof of Lemma \ref{sing}]
Fix $w$.  We solve \eqref{back} by writing an Ansatz 
$$ v = v^{\langle w \rangle} = v^{[w]} + \phi.$$

It is easily verified that $\phi$ will solve the difference equation

$$
i \phi_s + \phi_{yy} = -v^{[w]}_{yy} + s^{-1}F(\phi) $$
where $F = F_w$ denotes the quantity
$$
F(\phi) := (|v^{[w]} + \phi|^2 (v^{[w]} + \phi) - |v^{[w]}|^2 v^{[w]}).$$
We proceed by solving
the forward Cauchy problem with data $\phi(y,0) = 0$, 
in the sense that $\|\phi(s)\|_{\H^{K,K}} \to 0$ as $s \to 0^+$,  
rather than by specifying $v_1$ and solving the backwards Cauchy problem.
We then define $v_1(y)=v^{\langle w \rangle}(y,1)$.

We can write the difference equation in integral form as
\begin{equation}\label{diff}
\phi(s) = -\int_0^s U(s-s') v^{[w]}_{yy}(s')\ ds' +\int_0^s U(s-s')  
(s')^{-1}F(\phi(s'))\ ds'
\end{equation}
where $U(s) := \exp(-is \partial_{yy})$ is the free Schr\"odinger evolution operator.  We can solve this equation by setting up an iteration scheme
$$ \phi^{(k+1)}(s) = -\int_0^s U(s-s') v^{[w]}_{yy}(s')\ ds' +
\int_0^s U(s-s')  (s')^{-1}F(\phi^{(k)}(s'))\ ds'$$
with $\phi^{(0)} := 0$.

We claim inductively that
\be{phik-bound}
\| \phi^{(k)}(s) \|_{\H^{j,j}} \leq C \eps s (1 + |\log s|)^{10Kj + 10K}
\end{equation}
for all $0 < s \leq 1$, $0 \leq j \leq K$ and all $k$, where the constant $C$ is independent of $k$.  This is trivial for $k=0$.  Now assume it is proven for $k$.  To prove it for $k+1$, we observe from Minkowski's inequality and Lemma \ref{evolve} that
$$
\| \phi^{(k+1)}(s) \|_{\H^{j,j}}
\lesssim \int_0^s \| v^{[w]}_{yy}(s') \|_{\H^{j,j}}\ ds' +
 \int_0^s (s')^{-1}\| F(\phi^{(k)}(s')) \|_{\H^{j,j}}\ ds'.$$

By hypothesis we have $\|w\|_{\H^{K+2,K+2}} \lesssim \eps \ll 1$.  From \eqref{explicit}, \eqref{a-1}, \eqref{a-3} and the chain rule we observe the estimate
$$ \| v^{[w]} (s') \|_{\H^{K+2,K+2}} \lesssim \eps (1 + |\log(s')|)^{K+2}$$
and in particular that
$$ \| v^{[w]}_{yy}(s') \|_{\H^{K,K}} \lesssim \eps (1 + |\log(s')|)^{K+2}.$$
On the other hand, the $L^\infty$ norm of $v^{[w]}$ satisfies a bound
free of logarithms:
$$ \| v^{[w]} (s') \|_{L^\infty} \lesssim \eps.$$

By expanding out $F$ and using the Sobolev embedding $\| u \|_\infty \lesssim \| u \|_{\H^{1,1}}$ together with the bound on $\| v^{[w]} (s') \|_{L^\infty}$, we obtain
$$\| F(\phi^{(k)}(s')) \|_{L^2} \lesssim \| \phi^{(k)}(s') \|_{L^2} (\eps +  \| \phi^{(k)}(s') \|_{\H^{K,K}})^2.$$
If we expand out $F$ and use \eqref{a-1} (for multiplying $\phi^{(k)}$ with itself) and \eqref{a-2} (for multiplying anything with $v^{[w]}$), we obtain 
$$\| F(\phi^{(k)}(s')) \|_{\H^{j,j}} \lesssim (\| \phi^{(k)}(s')
\|_{\H^{j,j}} + \| \phi^{(k)}(s') \|_{\H^{j-1, j-1 }} (1 + |\log(s')|)^{K+2}) (\eps +  \| \phi^{(k)}(s') \|_{\H^{K,K}})^2
$$
for $1 \leq j \leq K$.  Inserting these bounds into the previous and applying the induction hypothesis we obtain
\bas
 \| \phi^{(k+1)}(s) \|_{\H^{j,j}} \lesssim & \eps s (1 + |\log(s)|)^{K+2}
+ C^3 \eps^3 s (1 + |\log(s)|)^{10Kj+10K} \\
& + C^3 \eps^3 s (1 + |\log(s)|)^{10K(j-1)+10K} (1 + |\log(s)|)^{K+2}.
\end{align*}
If $C$ is sufficiently large, and $\eps$ is sufficiently small depending on $C$, we may thus close the induction and obtain the desired bounds \eqref{phik-bound}. 

A standard variant of the above argument in fact shows that the iterates $\phi^{(k)}$ converge in $L^\infty_t((0,1]; \H^{K,K})$ to a solution $\phi$ to \eqref{diff} such that $\| \phi(s) \|_{\H^{K,K}} \to 0$ as $s \to 0^+$; indeed from \eqref{phik-bound} with $j=K$ we obtain \eqref{asym}. $v_1(y):=v^{\langle w\rangle}(y,1)$
then belongs to $\H^{K,K}$. 
By further standard arguments we can obtain the Lipschitz bound \eqref{lip}.
\end{proof}

\subsection{Decoherence and lack of scattering in $L^2$}

Now let $0 < \eps \ll 1$, and let $w$ be a non-zero $\H^{7,7}_y$ function with norm $O(\eps)$; for concreteness, let us take $w(y) := \eps e^{-y^2}$.  For any real number $a$ in the interval $[1/2,2]$ we apply Lemma \ref{sing} with $K=5$, and consider the function $v^{\langle aw \rangle}$.  At time $s=1$ this function depends continuously on $a$ in $\H^{5,5}_y$ norm.  However, we have the following {\it{decoherence property}} as $s \to 0^+$.

\begin{lemma}\label{eps}  If $|a|, |a'| = O(1)$ and $a \neq a'$ then
$$ \limsup_{s \to 0^+} \| v^{[aw]}(s) - v^{[a'w]}(s) \|_{L^2_y} \gtrsim 
(|a|+|a'|)\|w\|_{L^2}.$$
\end{lemma}

\begin{proof}
The conclusion is apparent\footnote{Indeed, it is immediate for each fixed $y$ that the integrand is of the order of $(|a|+|a'|)^2 |w(y)|^2$ on the average, and the claim follows by integrating in $s$ and using Fubini's theorem.} from 
\[
\|v^{[aw]}(s) - v^{[a'w]}(s) \|_{L^2}^2
= \int_{\R} 
|aw(y) e^{-i a^2|w(y)|^2 \log(s)} 
- a'w(y) e^{-i (a')^2|w(y)|^2 \log(s)}|^2\,dy.
\]
\end{proof}

From Lemma \ref{eps} and \eqref{asym} we thus have
\be{v-decohere}
\limsup_{s \to 0^+} \| v^{\langle aw \rangle } - v^{\langle a'w \rangle } \|_{L^2_y} \gtrsim 1
\end{equation}
for any $a \neq a'$, where the implicit constant depends on $w,a,a'$.  
This shows that the backwards Cauchy problem \eqref{back} is not uniformly well-posed in $L^2$ on the backwards interval $0 < s \leq 1$, even for initial data $v^{\langle w \rangle}$ with arbitrarily small $\H^{5,5}$ norm.  

The functions $v^{\langle aw \rangle}(s,y)$ solve the equation \eqref{back}.  By using the pseudo-conformal transformation \eqref{pc} we may thus construct solutions $u^{\langle aw \rangle } := \pc^{-1}(v^{\langle aw \rangle})$ to \eqref{nls}.  

\begin{remark}
To avoid confusion, we emphasize that $u^{\langle aw \rangle}$ is \emph{not} the solution to the Cauchy problem \eqref{nls} with initial datum $u(0) = aw$.
The ``initial condition'' involving $w$ is not posed at time $t=0$, but rather at time $t=\infty$, being given by the pullback of \eqref{asym} under the pseudo-conformal transformation.  However, since $u^{[aw]}(0) = \pc^{-1}(v^{[aw]})(1) = aw$, we expect from \eqref{asym} that $u(0)$ is in some sense ``close'' to $aw$.  Also, we do not know whether the map $w \to u^{\langle w \rangle}(0)$ is onto, even for Schwartz data; there may exist global solutions to \eqref{nls} whose asymptotic development does not resemble the one given here.  Certainly in the focusing case, the soliton solutions do not behave like the $u^{\langle aw \rangle}$.
\end{remark}

Together, Lemmas \ref{sing} and \ref{pc-convert} imply that
\be{uaw-bound}
\sup_{0 \leq t < \infty} \| u^{\langle aw \rangle}(t) \|_{H^5_x} \lesssim \eps.
\end{equation}
Furthermore, \eqref{asym} and Lemma \ref{pc-convert} imply
\begin{equation}
\label{uclose}
 \| u^{\langle aw \rangle}(t) - u^{[aw]}(t) \|_{H^5_x} 
\lesssim \eps (1+t)^{-1} \log^C (2 + t)
\end{equation}
for $0 \leq t < \infty$, where $u^{[aw]} = \pc^{-1} v^{[aw]}$ can be written explicitly as
\begin{equation}  \label{u[a]}
u^{[aw]}(t,x) = (1+t)^{-1/2} \exp(-ix^2/4(t+1)) a w(x/t) \exp(i a^2 |w(x/t)|^2 \log(1+t)).
\end{equation}
This implies that there is no scattering in $L^2$, or more precisely that 
that each element $u^{\langle aw\rangle}(t,x)$ of a large class of solutions 
fails to be asymptotically equal in $L^2_x$ norm to some solution of the
free Schr\"odinger evolution; this failure follows from the corresponding
failure for $u^{[aw]}(t,x)$ together with the inequality \eqref{uclose}. 
(Indeed, free $L^2$ solutions transform under $\pc$ to functions $v$ which can be continuously extended in $L^2$ to the time $s=0$, whereas the function $v^{[aw]}$ cannot be).  It may be that this failure of scattering can be repaired by modifying the free evolution appropriately, see \cite{ozawa}.

Together,
\eqref{lip} and Lemma \ref{pc-convert} give us
\be{soft-lip}
\| u^{\langle aw \rangle}(0) - u^{\langle a'w \rangle}(0) \|_{H^5_x} \lesssim \| v^{\langle aw \rangle }(1) - v^{\langle a'w \rangle }(1) \|_{\H^{5,5}_x} \lesssim \eps |a-a'|,
\end{equation}
while Lemma \ref{pc-convert} and \eqref{v-decohere} imply
\be{u-decohere}
\limsup_{t \to +\infty} \| u^{\langle aw \rangle }(t) - u^{\langle a'w \rangle }(t) \|_{L^2_x}
=  \limsup_{s \to 0^+} \| v^{\langle aw \rangle } - v^{\langle a'w \rangle } \|_{L^2_y} \gtrsim 1.
\end{equation}
We conclude in particular that the solution map $S$ for \eqref{nls} fails to be uniformly continuous from $L^2_x$ to $L^\infty_t((-\infty,\infty); L^2_x)$.  Thus the global well-posedness in $L^2$ in Theorem \ref{tsu} is not uniform in time.  
(Indeed, this argument shows one does not have uniformity for any $H^s$ for $s \geq 0$).

\section{Ill-posedness of NLS}\label{nls-sec}

We now indicate how the solutions $u^{\langle w \rangle}$ to \eqref{nls} constructed in the previous section disprove uniform continuity of the solution operator for defocusing NLS in $H^s$ for $s < 0$.  
Let $0 < \delta \ll \eps \ll 1$ and $T > 0$ be arbitrary.
We shall find two solutions $u = \phi^{\langle a \rangle}, \phi^{\langle a' \rangle}$ to \eqref{nls} such that
\begin{align}
\| \phi^{\langle a \rangle}(0) \|_{H^s_x}, \| \phi^{\langle a' \rangle}(0) \|_{H^s_x} &\lesssim \eps \label{0-small}\\
\| \phi^{\langle a \rangle}(0) - \phi^{\langle a' \rangle}(0) \|_{H^s_x} &\lesssim \delta \label{0-close} \\
\sup_{0 \leq t < T} \| \phi^{\langle a \rangle}(t) - \phi^{\langle a' \rangle }(t) \|_{H^s_x} &\gtrsim \eps.
\label{T-far}
\end{align}
This implies the solution map $S$ is not uniformly continuous from the ball $\{ u_0 \in H^s_x: \| u_0 \|_{H^s_x} \lesssim \eps \}$ to $L^\infty_t([0,T]; H^s_x)$, thus proving Theorem \ref{defocus-nls}.
It turns out that the subcritical case $-1/2 < s < 0$, the critical case $s=-1/2$, and the supercritical case $s < -1/2$ must be treated separately.

\subsection{The subcritical case}

Fix $s\in (-\tfrac12,0)$.
Let $N \gg 1$, $\lambda > 0$ be parameters to be chosen later, and let $K \geq 5$ be a large integer.  
Consider the functions
$$ \phi^{\langle a \rangle}(t,x) := 
\lambda e^{iN x } e^{i N^2 t} u^{\langle aw \rangle }
(\lambda^2 t, \lambda(x + 2tN))$$
where $a \in [1/2,2]$ is a parameter to be chosen later and $w$ is as in the previous section, with $\|w\|_{H^{K+2,K+2}}=O(\eps)$.  
We similarly define 
$\phi^{\langle a' \rangle}$ for some $a' \neq a$ also in $[1/2,2]$.

The Galilean and scale invariances of \eqref{nls} imply that $\phi^{\langle a \rangle}$ and $\phi^{\langle a' \rangle}$ are solutions of NLS.
Moreover, Lemma \ref{modulate} gives
$$ \|\phi^{\langle a \rangle}(0)\|_{H^s_x} 
\lesssim N^s \lambda^{1/2} \| u^{\langle aw \rangle}(0) \|_{H^K_x},$$
provided that its hypothesis $1\ll \lambda^{-1} N^{1+(s/K)}$ is
satisfied for all sufficiently large $N$.
We thus set $\lambda := N^{-2s}\gg 1$, so that \eqref{0-small} holds.
The condition $1\ll \lambda^{-1} N^{1+(s/K)}$
then becomes $1\ll N^{1+2s+(s/K)}$; hence for any $s>-1/2$ this condition will be obeyed for $N$ and $K$ sufficiently large. 

Similarly, from Lemma \ref{modulate} and \eqref{soft-lip} we have
$$ \|\phi^{\langle a \rangle}(0) - \phi^{\langle a' \rangle}(0) \|_{H^s_x} \lesssim \eps |a-a'|.$$
Thus \eqref{0-close} holds if $a$ and $a'$ are sufficiently close depending on $\delta$ (but still unequal).

By \eqref{u-decohere} there exists a time $t_0>0$ depending on $a, a'$ (but not on $N$, $\lambda$) such that
$$ \| u^{\langle aw \rangle}(t_0) - u^{\langle a'w \rangle}(t_0) \|_{L^2_x} \gtrsim \eps.$$
Fix this $t_0$.  From \eqref{uaw-bound} we have
$$ \| u^{\langle aw \rangle}(t_0)\|_{H^K_x}, \| u^{\langle a'w \rangle }(t_0) \|_{H^K_x} \lesssim \eps.$$
By Lemma \ref{modulate} we thus have
$$ \|\phi^{\langle a \rangle}(t_0 \lambda^{-2}) 
- \phi^{\langle a' \rangle}(t_0 \lambda^{-2})\|_{H^s_x} 
\sim N^s \lambda^{1/2} \eps 
= \eps.$$
If we choose $N$ large enough, we can make $t_0 \lambda^{-2} < T$, and so \eqref{T-far} follows.  This concludes the proof of Theorem \ref{defocus-nls} when $-1/2 < s < 0$.

\subsection{The critical case}

As $s \to -1/2$ in the subcritical argument above we see that we need the Galilean invariance less and less, and rely more on scaling.  Thus in the critical case $s = -1/2$ we expect to obtain ill-posedness purely by scaling the solutions $u^{\langle aw \rangle}$ used earlier.  However there is a slight difficulty in that we need some vanishing of the Fourier transform at the origin to make the $\dot H^{-1/2}$ norm converge.  Fortunately this can be easily achieved by making all the solutions odd.

We turn to the details.  We fix $K=5$ and $0 < \delta \ll \eps \ll 1$, $T > 0$ as before, and let $w$ be a function obeying $\|w\|_{\H^{K+2,K+2}} \lesssim \eps$.  We also assume that $w$ is \emph{odd} (in order to create vanishing at the frequency origin).  Then for any $a = O(1)$, the function $v^{[aw]}$ defined by \eqref{explicit} is also odd (in space), and an inspection of the argument in Lemma \ref{sing} shows that $v^{\langle aw \rangle}$ is similarly odd.  Inverting the pseudo-conformal transformation using Lemma \ref{sing} and Lemma \ref{pc-convert}, we see in particular that the initial datum
$u^{\langle aw \rangle}(0)$ is odd and obeys the estimates
$$
\| u^{\langle aw \rangle}(0) \|_{\H^{5,5}} \lesssim \eps.
$$
Taking Fourier transforms, we see in particular that
$$
|\partial_\xi^a \widehat{ u^{\langle aw \rangle}(0) }(\xi)| \lesssim (1+|\xi|)^{-1}$$
for $a = 0,1$.  Since $u^{\langle aw \rangle}(0)$ is odd, its Fourier transform vanishes at the origin, and so we thus have
$$
|\widehat{ u^{\langle aw \rangle}(0) }(\xi)| \lesssim \min(|\xi|, |\xi|^{-1}).$$
A similar argument yields
$$
|\widehat{ u^{\langle a'w \rangle}(0) }(\xi)| \lesssim |a-a'| \min(|\xi|, |\xi|^{-1}).$$
Now we pick a large parameter $\lambda \gg 1$ and redefine the functions $\phi^{\langle a \rangle}$ by
$$ \phi^{\langle a \rangle}(t,x) := \lambda u^{\langle aw \rangle}(\lambda^2 t, \lambda x)$$
(i.e.\ we just scale, and perform no additional Galilean transformation).
The estimates \eqref{0-small}, \eqref{0-close} then follow directly from the above pointwise Fourier transform estimates, while the proof of \eqref{T-far} is identical to the subcritical case.

\subsection{The supercritical case}

Finally, we address the supercritical case $s < - 1/2.$ We will once again
construct solutions $\phi^{\langle aw \rangle}, \phi^{\langle a' w \rangle}$ of
\eqref{nls} satisfying the conditions \eqref{0-small}, \eqref{0-close}
and \eqref{T-far}. The construction in this case will rely on an approximation
by the zero dispersion limit of NLS rather than the modified scattering
solutions above. The method used here is rather general and might potentially
be used to establish ill-posedness in the supercritical regime for many equations.

We begin by considering the {\it{small dispersion defocusing NLS}}
initial value problem
\begin{equation}
  \label{NLSdelta}
 \left\{
\begin{aligned}
-i v_t+ \delta^2 v_{xx} &= |v|^2 v \\
v(0, x) &= f(x)
\end{aligned}
\right.
\end{equation}
where the {\it{dispersion parameter}} $\delta$ satisfies $0 < \delta\leq 1$.
We assume that the initial datum $f$ belongs to the Schwartz space
${\mathcal S}$.
In the limit $\delta=0$, a solution is $v(t,x) = f(x)e^{it|f(x)|^2}$.

We will consider solutions of \eqref{NLSdelta} with initial data
$aw(x)$, where $a\in\complex$ varies freely within the unit ball,
and $w\in{\mathcal S}$ is fixed.  Set $v^{[aw]}(t,x) = aw(x)e^{it|aw(x)|^2}$.

\begin{lemma} \label{smalldispersionapprox}
Let $w\in{\mathcal S}$ and $N\in\naturals$ be given. Then there exist
constants $C,q\in\R^+$ and
a lifespan function $T:(0,1]\mapsto\R^+$, such that $T(\delta)\to+\infty$
as $\delta\to 0$, with the following property.
For each $\delta\in(0,1]$
and each $a \in\complex$ satisfying $|a|\le 1$,
there exists a solution $v = v^{\langle aw,\delta\rangle}(t,x)$
in $C^0(\R, H^{N,N})$
of the small dispersion NLS initial value problem
\eqref{NLSdelta} with dispersion coefficient $\delta^2$
satisfying
\begin{equation}
\|v^{\langle aw,\delta\rangle}(t,\cdot)
- v^{[aw]}(t,\cdot)\|_{H^{N,N}(\R)}
\le C\delta^q
\end{equation}
uniformly for all $|t|\le T(\delta)$.
\end{lemma}

The proof is based on the standard energy method, very much like
the proofs of Lemma~\ref{sing} and Theorem~\ref{perturb}. Existence
of a solution for all time is well-known, since a simple change
of variables reduces matters to the case $\delta=1$. Plugging
the initial approximation $v^{[aw]}$ into the differential equation,
one finds that $-iv^{[aw]}_t+\delta^2 v^{[aw]}_{xx}
- |v^{[aw]}|^2 v^{[aw]}$ is $O(\delta^2)$ in $C^0(\R,H^{N,N})$ norm.
Seeking a solution of the form $v^{[aw]}+u$,
one then analyzes $\partial_t\|u(t\cdot) \|_{H^{k,k}}^2$.
Only {\em upper} bounds for $\delta$, together with the condition
$\delta^2\in\R$, are required to carry this out. The details of the proof
are left to the reader.

From the remainder bound of Lemma~\ref{smalldispersionapprox}
together with an elementary
comparison of the approximate solutions $v^{[aw]}$ there follows
a decoherence property: For any distinct $a,a'\in\complex$
satisfying $|a|,|a'|\le 1$ and any $r\in\R$, there exists $\eta>0$
such that for each $\delta\le\eta$,
there exists $\tau>0$ satisfying $\tau\le C\big||a|^2-|a'|^2 \big|^{-1}$
for which
\begin{equation}
\|v^{\langle aw,\delta\rangle}(\tau,\cdot)
- v^{\langle a'w,\delta\rangle}(\tau,\cdot)\|_{H^r(\R)}
\ge c|a|+c|a'|.
\end{equation}

Consider next the functions
\begin{align*}
g^{\langle aw,\delta\rangle}(t,x) &= v^{\langle aw,\delta\rangle}(t,\delta x),
\\
g^{[aw]}(t,x) & = v^{[aw]}(t,\delta x) = aw(\delta x) e^{i|aw(\delta x)|^2}.
\end{align*}
$g^{\langle aw,\delta\rangle}$ is an exact solution of the
defocusing NLS equation \eqref{nls}, with initial datum $aw(\delta x)$.

As they stand,
these solutions are unsuitable for an ill-posedness argument, because
the $H^s$ norm of the initial datum tends to infinity as $\delta\to 0$.
However, further solutions may obtained via the scaling symmetry
of the equation: For each $\lambda\in\R^+$,
\begin{equation}
\phi^{\langle aw,\delta,\lambda\rangle}(t,x) =
\lambda\rp g^{\langle aw,\delta\rangle}(\lambda^{-2}t,\lambda\rp x)
\end{equation}
is also a solution of \eqref{nls}, with initial datum $\lambda\rp a
w(\lambda\rp\delta x)$.
For $\phi^{\langle a,\delta\rangle}$ we have the approximation
$\phi^{[a,\delta]}(t,x) = \lambda\rp a
w(\lambda\rp\delta x) e^{i\lambda^{-2}t|aw(\delta x)|^2}$.

Suppose now that $s<-1/2$, so that the Sobolev space $H^s$ is
supercritical for the NLS equation.
We wish to choose $\lambda$ as a function of $\delta$, so that
the $H^s$ norm of the initial datum $\lambda\rp a
w(\lambda\rp\delta x)$ is $\sim|a|$, uniformly as $\delta\to 0$.
A simple calculation shows that in the {\em homogeneous}
Sobolev space $\dot H^s$, such a normalization is achieved by
taking
\begin{equation} 
\lambda = \delta^\gamma \ \ \ \ \text{where}\ \gamma = \gamma(s)
= \frac{-2s+1}{-2s-1}.
\end{equation}
For all $s<-1/2$, $\gamma>1$.

Before proceeding, we impose a restriction on the Schwartz class function
$w$ which will be needed below:
\begin{equation}
\hat w(\xi) = O(|\xi|^\kappa)\ \text{as}\ \xi\to 0,
\end{equation}
where $\kappa$ is a large positive integer, depending on $s$,
to be specified below.
Of course, we also require that $w$ not vanish identically;
more specifically, we require
that $\hat w$ not vanish identically on the interval $I=[1,2]$.

We now define a two-parameter family of solutions of the NLS equation
by 
\[ \phi^{\langle a,\delta\rangle}(t,x)
= \phi^{\langle aw,\delta,\lambda\rangle}(t,x) 
\ \ \ \text{where}\  \lambda = \delta^{\gamma(s)}.
\] 
The corresponding initial data are $\varphi_{a,\delta}(x)
= \delta^{-\gamma}a w(\delta^{1-\gamma}x)$.
Since $\delta^{1-\gamma}\to\infty$ as $\delta\to 0$, 
and $\hat w(0)=0$, these initial data are composed primarily
of higher-frequency Fourier modes as $\delta$ becomes smaller.
The corresponding explicit approximations are
$\phi^{[a,\delta]}(t,x) = \delta^{-\gamma}a w(\delta^{1-\gamma}x)
e^{i\delta^{-2\gamma} t |aw(\delta^{1-\gamma}x)|^2}$.
The functions $\phi^{\langle a,\delta\rangle}(t,x)$ satisfy
the NLS equation globally in time; however, we have reasonably
good control over them only for short times
$|t|\le \delta^{2\gamma}T(\delta)$.

We next verify that these initial data are approximately normalized
in the inhomogeneous Sobolev spaces $H^s$.
Since $\widehat{\varphi_{a,\delta}}(\xi)
= a \delta\rp \hat w(\delta^{\gamma-1}\xi)$,
\begin{equation*}
\|\varphi_{a,\delta}\|_{H^s}^2
=
|a|^2
\delta^{-1-\gamma} \int_{\R} \big|\hat w(\xi) \big|^2
(1 + |\delta^{1-\gamma}\xi|^2)^s\,d\xi.
\end{equation*}
Split the region of integration into two parts.
The contribution of the region $|\xi|\ge \delta^{\gamma-1}$
is
\[
\sim |a|^2 \delta^{-1-\gamma}\delta^{2s(1-\gamma)}
\int_{|\xi|\ge \delta^{\gamma-1}} |\hat w(\xi)|^2\,d\xi.
\]
Since $(1+\gamma)/(1-\gamma) = 2s$, and since $\delta^{\gamma-1}
\to 0$ as $\delta\to 0$, this is
\[
\sim c|a|^2 \|w\|_{L^2}^2\ \text{ as }\ \delta\to 0.
\]
The contribution of the region $|\xi|\le \delta^{\gamma-1}$ is
\[
\le C |a|^2 \delta^{-1-\gamma} \delta^{(\gamma-1)(1+ 2 \kappa)}.
\]
Since $\gamma>1$, for any given $s<-1/2$ there exists $\kappa$
such that this last expression is $O(|a|^2 \delta)$ as $\delta\to 0$.
Thus we conclude that
\[
\|\varphi_{a,\delta}\|_{H^s}^2 = c_0|a|^2\|w\|_{L^2}^2
+ O(|a|^2\delta) \ \text{as}\  \delta\to 0.\]

The next step is to argue that $\phi^{\langle a,\delta\rangle}$
is nearly equal to $\phi^{[a,\delta]}$.
Here arises a complication, because there is no analogue for
positive time of the condition 
$\hat w(\xi) = O(|\xi|^\kappa)$.\footnote{We hope to exploit this
in a future paper to establish a more dramatic form of ill-posedness,
in the supercritical case.}
What does follow directly from Lemma~\ref{smalldispersionapprox}
is that
\begin{equation}  \label{minus}
\int_{|\eta|\ge 1}
\big| \widehat{\phi^{\langle a,\delta\rangle}}(t,\eta)
-
\widehat{\phi^{[a,\delta]}}(t,\eta)\big|^2
(1+|\eta|^2)^{s}\,d\eta
\le C\delta^{2q}
\end{equation}
uniformly for $|t|\le \delta^{2\gamma}T(\delta)$.
Indeed, in this region $(1+|\eta|^2)^s\sim |\eta|^{2s}$,
and a change of variables as in the above calculation of the norm of
$\varphi_{a,\delta}$, together with the bound from
 Lemma~\ref{smalldispersionapprox}, yields the result.
The factor of $\delta^{2\gamma}$ in the upper bound for $|t|$
arises because the definition of $\phi^{\langle a,\delta\rangle}$
involves a rescaling of time.
Here the Fourier transform is taken in the variable $x$, for each
time $t$.

Next,
\[
\int_{1< \delta^{\gamma-1}\eta<2}
\big| \widehat{\phi^{[a,\delta]}}(t,\eta)
-
\widehat{\phi^{[a',\delta]}}(t,\eta)\big|^2
(1+|\eta|^2)^{s}\,d\eta 
\sim \int_{1<\xi<2}
\big| 
(awe^{i\delta^{-2\gamma}t|aw|^2})^\wedge (\xi)
-
(a'we^{i\delta^{-2\gamma}t|a'w|^2})^\wedge(\xi)
\big|^2\,d\xi
\]
in the sense that each side is dominated by a universal constant
multiple of the other.
The right-hand side equals
\[
|a|^2
\int_{1<\xi<2}
\big|
(we^{i\delta^{-2\gamma}t|aw|^2})^\wedge(\xi)
-
(we^{i\delta^{-2\gamma}t|a'w|^2})^\wedge(\xi)
\big|^2\,d\xi
+O(|a-a'|),
\]
since $w\in{\mathcal S}$ and $|a|,|a'|\le 1$.

Set $F(s) = (we^{is|w|^2})^\wedge$,
regarded as an element of $L^2(I)$, where $I=[1,2]$.
$F$ is an entire holomorphic function, depending periodically
on the real part of $s$ with period $2\pi$.
By choosing a generic Schwartz function $w$
(still satisfying $\hat w(\xi) = O(|\xi|^\kappa)$)
we may ensure that $F$ is nonconstant.
From this it follows by elementary reasoning that
whenever $b\ne b'\in [0,1]$,
there exist $s\in\R$ satisfying $0<s\le C|b-b'|\rp$
and $b''\in[0,1]$ such that $\tfrac12\le|b-b''|/|b-b'|\le 1$,
such that $\|F(bs)-F(b''s)\|_{L^2(I)}\ge c_0>0$.
Hence whenever $|a|\ne|a'|$,
there exist $0<t^*\le C\delta^{2\gamma}\big||a|^2-|a'|^2\big|\rp$
and $a''$ such that $\tfrac12 \le|a-a''|/|a-a'|\le 1$
and
\[
\int_{1< \delta^{\gamma-1}\eta<2}
\big| \widehat{\phi^{[a,\delta]}}(t^*,\eta)
-
\widehat{\phi^{[a'',\delta]}}(t^*,\eta)\big|^2
(1+|\eta|^2)^{s}\,d\eta
=
c_0|a|^2 + O(|a-a''|),
\]
where $c_0$ is a nonvanishing constant depending only on $w$.

If in addition $\big||a|^2-|a''|^2\big|\rp\le cT(\delta)$
for a sufficiently small constant $c>0$, as is the case
for all sufficiently small $\delta$,
then we may combine this with \eqref{minus}
to conclude that there exists $0<t^*\le C\delta^{2\gamma}
\big||a|^2-|a''|^2 \big|\rp$ such that
\[
\int_{1< \delta^{\gamma-1}\eta<2}
\big| \widehat{\phi^{\langle a,\delta\rangle}}(t^*,\eta)
-
 \widehat{\phi^{\langle a'',\delta\rangle}}(t^*,\eta)
\big|^2 (1+|\eta|^2)^s\,d\eta
\ge c_0|a|^2 + O(|a-a''|) + O(|a|\delta^{2a}).
\]

Fix any $a,a'$ such that
$0<|a|<1$, $|a|\ne |a'|$ and $\big||a'|-|a|\big|/|a|$ is
less than a small fixed constant.
To these associate a substitute parameter $a''$, as above.
Then for every sufficiently small $\delta$
there exists $t^*=O(\delta^{2\gamma})$ such that
\[
\|{\phi^{\langle a,\delta\rangle}}(t^*,\cdot)
-
{\phi^{\langle a'',\delta\rangle}}(t^*,\cdot)\|_{H^s(\R)}
\ge \tfrac12 c_0 |a|^2. 
\] 
Because $\delta^{2\gamma}\to 0$ as $\delta\to 0$,
this means
that the solution operator fails to be uniformly continuous,
even when restricted to any small neighborhood of the origin in $H^s$.

\noindent 
{\em Remark.\/}
The essential feature of the supercritical case used in this proof
is that the same scaling transformations which reduce the
(homogeneous) $\dot H^s$ norm also {\em contract} the time variable.
~

\subsection{Remarks}

\begin{remark}
For this discussion we stay in the subcritical case $s > -1/2$ and allow implicit constants to depend on $\eps$ and $T$.  A more careful inspection of Lemma \ref{eps} shows that one begins to have decoherence at time $s \sim \exp(-C/|a-a'|)$.  In our application $|a-a'| \sim \delta$.  Chasing through all the constants we 
obtain $N \sim \exp(C/\delta)$.  Thus, our counterexample is quite weak in the sense that we need to go out to frequencies $\sim \exp(C/\delta)$ to obtain a failure of uniform continuity at uncertainty $\delta$.  In comparison, the soliton-based arguments in \cite{kpv:counter} only require that one go out to frequencies $\sim \delta^{-C}$ to achieve a similar result.
It would potentially be interesting if this weakness reflected a genuine
feature of the equation.
\end{remark}

\begin{remark} 
One may informally compare the results here and those in \cite{kpv:counter} from the perspective of complete integrability.  The NLS equation is completely integrable and can be studied by inverse scattering techniques.  For the focusing NLS equation, a general solution can be split into a ``multisoliton'' component, which eventually resolves into a collection of disjoint solitons, and a ``dispersive'' component, which eventually decays \cite{Segur:NLS}, \cite{ZM:NLS}.  In the defocusing equation there are far fewer solitons, and the behaviour is mostly dispersive.  The NLS results in \cite{kpv:counter} can be viewed as a statement that the soliton component of NLS is badly behaved in negative Sobolev spaces; the results here say (informally speaking) that the dispersive component is also badly behaved in these spaces (though to a lesser degree - see previous remark).  One can also view the results on KdV and mKdV in this way.  However we emphasize that our methods here do not require complete integrability or explicit travelling wave solutions, and should extend to other, non-integrable equations such as the nonlinear wave equation.
\end{remark}

\section{A review of local well-posedness for mKdV}\label{mkdv-lwp-sec}

We now turn our attention to the modified KdV (mKdV) equation \eqref{mkdv}.  We begin by reviewing the local well-posedness theory of \cite{kpv:gkdv} for \eqref{mkdv} at the endpoint regularity $H^{1/4}$. We recall from \cite{kpv:gkdv} the following linear estimates for the Airy equation.

\begin{theorem}\label{airy}\cite{kpv:gkdv}  Suppose that $u$ solves the inhomogeneous problem
$$ u_t + u_{xxx} = F; \quad u(0,x) = u_0(x)$$
on the slab $[0,T] \times \R$.  Then $u$ satisfies 
the smoothing estimate (cf. \cite{kato}, \cite{kruzkhov}, \cite{kpv:gkdv})
\be{kato}
 \| u_x \|_{L^\infty_x(L^2_t)} \lesssim \| u_0 \|_{L^2_x} + \int_0^T \| F(t) \|_{L^2_x}\ dt
\end{equation}
and the maximal function estimates
\begin{align}
\| u \|_{L^4_x(L^\infty_t)} &\lesssim \| u_0 \|_{H^{1/4}} + \int_0^T \| F(t) \|_{H^{1/4}}\ dt \label{max-14}\\
\| u \|_{L^2_x(L^\infty_t)} &\lesssim \| u_0 \|_{H^{3/4+\epsilon}} + \int_0^T \| F(t) \|_{H^{3/4+\epsilon}}\ dt. \label{max-34}
\end{align}
\end{theorem}

The estimate \eqref{max-34} is not really needed for the present discussion, but will be used
to deal with an mKdV-like system in Section \ref{general-sec}.

For any time interval $I = [t_0, t_0+T]$, let $X = X(I \times \R)$ denote the norm
\begin{equation} \label{Xnorm}
\| u \|_{X(I \times \R)} := \| u \|_{L^\infty_t ( H^{1/4}_x)} + 
\| u \|_{L^4_x(L^\infty_t)}  
+ \| D^{-1/2-\epsilon} u \|_{L^2_x(L^\infty_t)}
+
\| \partial_x D^{1/4} u \|_{L^\infty_x(L^2_t)} 
\end{equation}
on the spacetime slab $I \times \R$, where $D := \sqrt{1+(-\Delta)}$.  
From the above theorem and energy estimates there follows the inequality
\be{energy}
\| u \|_{X(I \times \R)} \lesssim \| u(t_0) \|_{H^{1/4}_x} + 
\int_{t_0}^{t_0+T} \| (\partial_t + \partial_{xxx}) u(t) \|_{H^{1/4}_x}\ dt.
\end{equation}

In \cite{kpv:gkdv} the following trilinear estimate was proven (see also the proof of Proposition \ref{iteration}):

\begin{theorem}\label{est}\cite{kpv:gkdv}  On any spacetime slab $I \times \R$, we have
$$ 
\| D^{1/4} (uv w_x) \|_{L^2_x(L^2_t)} \lesssim \|u\|_X \|v\|_X \|w\|_X 
$$
\end{theorem}

By combining this estimate with \eqref{energy} one can obtain local well-posedness for mKdV in $H^{1/4}$; see \cite{kpv:gkdv}.  One can also use these estimates in a standard manner to obtain the following perturbation result for the mKdV flow in $H^{1/4}$:

\begin{lemma}\label{perturbation}  Suppose that $u$ is a smooth solution to the mKdV equation \eqref{mkdv}, and suppose that $v$ is an approximate Schwartz solution to mKdV in the sense that
$$ v_t + v_{xxx} = 6 v^2 v_x - E$$
for some error function $E$.  Let $t_0$ be a time, and let $e$ be the solution to the inhomogeneous problem
$$ e_t + e_{xxx} = E; \quad e(t_0) = 0.$$
Suppose that we have the estimates
$$ \| u(t_0) \|_{H^{1/4}_x}, \| v(t_0) \|_{H^{1/4}_x} \lesssim \eps; \quad \| e \|_{X([t_0,t_0+1] \times \R)} \lesssim \eps$$
for some sufficiently small absolute constant $0 < \eps \ll 1$.  Then we have
$$ \| u-v \|_{X([t_0,t_0+1] \times \R)} \lesssim \| u(t_0)-v(t_0)\|_{H^{1/4}_x} + \| e \|_{X([t_0,t_0+1] \times \R)}.$$
In particular we have
\begin{equation} \label{pert-ineq}
 \sup_{t_0 \leq t \leq t_0+1} \| u(t) - v(t) \|_{H^{1/4}_x}
\lesssim \| u(t_0)-v(t_0)\|_{H^{1/4}_x} + \| e \|_{X([t_0,t_0+1] \times \R)}.
\end{equation}
\end{lemma}
In other words, any function $v$ which approximately satisfies mKdV in the above sense 
stays close to the exact mKdV flow.

\begin{proof} In this proof we work entirely on the spacetime slab $[t_0,t_0+1] \times \R$.
Write the equation for $v$ in integral form as
$$ v(t) = U(t-t_0) v(t_0) - e(t) + \int_{t_0}^t U(t'-t_0) (6 v^2 v_x )(t')\ dt'$$ 
where $U(t) := \exp(-t\partial_{xxx})$ is the free Airy evolution operator.  Taking $X$ norms of both sides and using \eqref{energy} we obtain
$$ \| v\|_X \lesssim \| v(t_0) \|_{H^{1/4}_x} + \| e \|_X + \| v^2 v_x \|_{L^1_t(H^{\frac{1}{4}}_x)}. $$ 
By H\"older's inequality
we may estimate the $L^1_t(H^{1/4}_x)$ norm by the $L^2_t(H^{1/4}_x)$ norm. 
Using Theorem \ref{est} we thus have
$$ \| v\|_X \lesssim \| v(t_0) \|_{H^{1/4}_x} + \| e \|_X + \| v \|_X^3.$$ 
If $\eps$ is sufficiently small, we thus deduce via a continuity argument that
\be{vx-eps}
\|v\|_X \lesssim \eps.
\end{equation}

If we write $u = v+w$, then $w$ satisfies the equation
$$ w_t + w_{xxx} = 2 \partial_x( 3wv^2 + 3w^2v + w^3 ) + E; \quad w(t_0) = u(t_0) - v(t_0)$$
which can be written in integral form as
$$ w(t) = U(t-t_0)(u(t_0)-v(t_0)) + e(t) + \int_{t_0}^t U(t'-t_0) 2 \partial_x( 3wv^2 + 3w^2v + w^3 )(t')\ dt'.$$
We again take $X$ norms and use \eqref{energy} to obtain
$$ \| w\|_X \lesssim \| u(t_0)-v(t_0) \|_{H^{1/4}_x} + \| e\|_X + \| D^{1/4} \partial_x(3wv^2 + 3w^2v + w^3) \|_{L^1_t([t_0,t_0+1]; L^2_x)}.$$
Again we estimate the $L^1_t(L^2_x)$ norm by the $L^2_t(L^2_x)$ norm and use Theorem \ref{est} to obtain
$$  \| w\|_X \lesssim \| u(t_0)-v(t_0) \|_{H^{1/4}_x} + \| e\|_X + \|w\|_X (\|w\|_X + \|v\|_X)^2.$$
By another continuity argument and \eqref{vx-eps} we obtain the desired result, if $\eps$ is sufficiently small.
\end{proof}

\section{A crude proof of Theorem \ref{defocus-mkdv}}\label{crude-sec}

In this section we give a proof of Theorem \ref{defocus-mkdv} relying on the rather crude perturbation result of Lemma \ref{perturbation}.  While this 
suffices to establish ill-posedness, it is quite poor quantitatively, 
and in the next two sections we shall give an argument which is similar in strength to the proof of Theorem \ref{defocus-nls}.

The first step is to construct $H^{1/4}$-normalized solutions whose asymptotic development can be controlled for relatively long periods of time.  Then, as in Section \ref{nls-sec},
 we shall use a scaling argument to demonstrate ill-posedness below $H^{1/4}$.

We recall from Section \ref{pseudo-sec} the global solutions $u^{\langle aw \rangle}$ to the NLS equation \eqref{nls} for all $a \in [1/2,2]$, where $w(x) = \eps \exp(-x^2)$ for some parameter $0 < \eps \ll 1$ to be chosen later.  As foreshadowed in the introduction, we shall use these NLS solutions to construct approximate solutions $V^{[a]}$ to \eqref{mkdv}, defined using the change of variables
$$ (s,y) := (t, \frac{x + 3N^2t}{3^{1/2} N^{1/2}})$$
by
\begin{equation} \label{V[a]defn}
V^{[a]}(t,x) := \sqrt{\frac{2}{3N}} \Re~  e^{iNx} e^{iN^3 t} u^{\langle aw \rangle}(s,y),
\end{equation}
where $N \gg 1$ is a large parameter to be chosen later.  From \eqref{uaw-bound} and Lemma \ref{modulate} we have
$$ \sup_{0 \leq t < \infty} \| V^{[a]}(t) \|_{H^{1/4}_x} \lesssim \eps.$$
Now we show that $V^{[a]}$ is an approximate solution to mKdV.  A straightforward computation shows that
$$ (\partial_t + \partial_{xxx}) V^{[a]}(t,x) = 
\sqrt{\frac{2}{3N}} \Re~  e^{iNx} e^{iN^3 t} 
(\partial_s + i\partial_{yy} + 3^{-3/2} N^{-3/2} \partial_{yyy}) 
u^{\langle aw \rangle}(s,y)$$
and that
\bas 2 \partial_x (V^{[a]}(t,x)^3) =&
\sqrt{\frac{2}{3N}} N^{-1} \partial_x 
\Big(\Re~ e^{iNx} e^{iN^3 t} |u^{\langle aw \rangle}(s,y)|^2 u^{\langle aw \rangle}(s,y) \\
&+ \frac{1}{3}
\Re~ e^{3iNx} e^{3iN^3 t} u^{\langle aw \rangle}(s,y)^3\Big) \\
=& 
\sqrt{\frac{2}{3N}} N^{-1} 
\Bigl[ \Re~ iN e^{iNx} e^{iN^3 t} |u^{\langle aw \rangle}(s,y)|^2 u^{\langle aw \rangle}(s,y) \\
&+ \Re~ e^{iNx} e^{iN^3 t} 3^{-1/2} N^{-1/2} 
\partial_y \big(|u^{\langle aw \rangle}(s,y)|^2 u^{\langle aw \rangle}(s,y)\big) \\
&+ 
\Re~ iN e^{3iNx} e^{3iN^3 t} u^{\langle aw \rangle}(s,y)^3 \\
&+ \frac{1}{3}
\Re~ e^{3iNx} e^{3iN^3 t} 3^{-1/2} N^{-1/2} \partial_y(u^{\langle aw \rangle}(s,y)^3) \Bigr].
\end{align*}
Since $u^{\langle aw \rangle}$ is a solution of \eqref{nls}, the main terms
of the preceding two equations agree, leaving
$$ (\partial_t + \partial_{xxx}) V^{[a]}(t,x) = 2 \partial_x (V^{[a]}(t,x))^3 - E$$
where the error term $E$ is a linear combination of the real and imaginary parts of the expressions 
\bas
E_1 &:= N^{-2} e^{iNx} e^{iN^3 t} u^{\langle aw \rangle}_{yyy}(s,y) \\
E_2 &:= N^{-2} e^{iNx} e^{iN^3 t} \partial_y(|u^{\langle aw \rangle}(s,y)|^2 u^{\langle aw \rangle}(s,y)) \\
E_3 &:= N^{-2} e^{3iNx} e^{3iN^3 t} \partial_y((u^{\langle aw \rangle}(s,y))^3) \\
E_4 &:= N^{-1/2} e^{3iNx} e^{3iN^3 t} (u^{\langle aw \rangle}(s,y))^3.
\end{align*}

We now give some estimates for $E_1, \ldots, E_4$, or more precisely for the solution to the inhomogeneous problem with these forcing terms.

\begin{lemma}\label{t-bone}  Let $t_0 \geq 0$. For each $j\in\{1,2,3,4\}$  
let $e_j$ be the unique solution to the problem
$$ (\partial_t + \partial_{xxx}) e_j = E_j; \quad e_j(t_0) = 0.$$
Then 
$$ \| e_j \|_{X([t_0, t_0+1] \times \R)} \lesssim \eps N^{-3/2}$$
where $X([t_0, t_0+1] \times \R)$ norm is as defined in \eqref{Xnorm}.
The inequality holds uniformly in $t_0$.
\end{lemma}

\begin{proof}
First suppose that $j\in\{1,2,3\}$.  
Then by \eqref{energy} it will suffice to show that
$$ \sup_t \| E_j(t) \|_{H^{1/4}_x} \lesssim \eps N^{-3/2}.$$
However, from \eqref{uaw-bound} and the fact that $H^k_y$ is closed under multiplication for all $k \geq 1$, we see that the functions $u^{\langle aw \rangle}_{yyy}$, $\partial_y(|u^{\langle aw \rangle}|^2 u^{\langle aw \rangle})$, and $\partial_y((u^{\langle aw \rangle})^3)$ all have an $H^1_y$ norm of $O(\eps)$.  The above claim then follows 
directly from Lemma \ref{modulate}. 

These arguments do not work for $E_4$ as this term does not contain enough negative powers of $N$, and one would only obtain a bound such as $O(\eps^3)$ which is insufficient for our purposes.   To do better we take advantage of the oscillation $e^{3iNx} e^{3iN^3 t}$, using the fact that the frequency $(\tau,\xi) = (3N^3,3N)$ is quite far away from the cubic $\tau = \xi^3$.  

A computation shows that
$$ (\partial_t + \partial_{xxx}) E_4 = -24 i N^3 E_4 + f$$
where $f$ is a linear combination of
\bas
f_1 &:= N^{-1/2} e^{3iNx} e^{3iN^3 t} \partial_s (u^{\langle aw \rangle}(s,y)^3) \\
f_2 &:= N e^{3iNx} e^{3iN^3 t} \partial_y (u^{\langle aw \rangle}(s,y)^3) \\
f_3 &:= N^{-1/2} e^{3iNx} e^{3iN^3 t} \partial_{yy} (u^{\langle aw \rangle}(s,y)^3) \\
f_4 &:= N^{-2} e^{3iNx} e^{3iN^3 t} \partial_{yyy} (u^{\langle aw \rangle}(s,y)^3).
\end{align*}
Rewriting the above as
$$  (\partial_t + \partial_{xxx}) (e_4 - \frac{1}{-24iN^3} E_4) = - \frac{1}{-24iN^3} f$$
and using \eqref{energy} we thus have
$$ \| e_4 \|_{X} \lesssim
N^{-3} ( \| E_4 \|_X + \| E_4(t_0) \|_{H^{1/4}_x} + \sum_{j=1}^4 \int_{t_0}^{t_0+1} \| f_j (t)\|_{H^{1/4}_x}\ dt ).$$
To estimate the functions $f_j$, we again observe from \eqref{uaw-bound} that the functions $\partial_s (u^{\langle aw \rangle}(s,y)^3)$ and $\partial^m_y (u^{\langle aw \rangle}(s,y)^3)$ are in $H^1_y$ for $m=1,2,3$ (for the $\partial_s$ derivative, we use the Leibnitz rule followed by \eqref{nls} to convert it to spatial derivatives), 
and so these terms are $O(N^{+3/2})$, as a consequence of Lemma \ref{modulate}. 
Taking the factor of $N^{-3}$ which multiplies each $f_j$ into account, we
obtain the desired bound $O(N^{-3/2})$.

It remains to treat the $E_4$ terms; 
since $\|E_4\|_X$ controls $\|E_4(t_0) \|_{H^{1/4}_x}$ it will suffice to show that
$$ \| E_4 \|_X \lesssim \eps N^{3/2}.$$
This can be done by direct computation (possibly using Sobolev embedding to first replace the mixed spacetime norms by unmixed norms) but one can also exploit (a modulated version of) \eqref{energy}.  Define $\tilde E_4$ by
$$ \tilde E_4 := e^{24iN^3t} E_4 = N^{-1/2} e^{3iNx} e^{27iN^3 t} u^{\langle aw \rangle}(s,y)^3;$$
observe that $\tilde E_4$ has the same $X$ norm as $E_4$, but unlike $E_4$, the function $\tilde E_4$ lives near the cubic $\tau = \xi^3$ in frequency space.  By \eqref{energy} it will suffice to show that
$$ \|\tilde E_4(t_0) \|_{H^{1/4}_x} + \int_{t_0}^{t_0+1} \| (\partial_t + \partial_{xxx}) \tilde E_4(t) \|_{H^{1/4}_x}\ dt \lesssim \eps N^{3/2}.$$
The first term is easily checked by Lemma \ref{modulate} and \eqref{uaw-bound}, so we turn to the latter.  We can expand $(\partial_t + \partial_{xxx}) \tilde E_4(t)$ as a linear combination of $e^{24iN^3 t} f_j$ for $j=1,2,3,4$.  But these terms have already been shown to be $O(N^{+3/2})$  in $L^\infty_t(H^{1/4}_x)$, and we are done.
\end{proof}

Let $V^{\langle a \rangle}$ be the global smooth solution to \eqref{mkdv} with initial 
datum $V^{\langle a \rangle}(0) = V^{[a]}(0)$.  Lemma \ref{t-bone}, Lemma \ref{perturbation}, and an easy induction argument give
\begin{equation}  \label{Vdifference}
\| V^{\langle a \rangle}(t) - V^{[a]}(t) \|_{H^{1/4}_x}  \lesssim \eps C^t N^{-1}
\end{equation}
for all $0 \leq t \ll \log N$.  
Indeed, applying the two lemmas gives
\begin{align*}
\sup_{j\le t\le j+1} \| V^{\langle a \rangle}(t) - V^{[a]}(t) \|_{H^{1/4}_x} 
&\le C
\| V^{\langle a \rangle}(j) - V^{[a]}(j) \|_{H^{1/4}_x} 
+ 
C\eps N^{-1}
\\
& \le C
\sup_{j-1\le t\le j} \| V^{\langle a \rangle}(t) - V^{[a]}(t) \|_{H^{1/4}_x} 
+ C\eps N^{-1}
\end{align*}
and the desired conclusion follows by induction on $j$
for $j\lesssim \log(N)$;
for $j=0$ one uses instead the identity
$V^{\langle a \rangle}(0) = V^{[a]}(0)$.

We conclude that for any $\eta>0$,
there exists a constant $c_0>0$ such that
\be{09}
\sup_{0 \leq t \le c_0\log N} 
\| V^{\langle a \rangle}(t) - V^{[a]}(t) \|_{H^{1/4}_x}  
\lesssim \eps N^{-1+\eta},
\end{equation}
uniformly for all $N\ge 2$.

We have thus constructed a one-parameter family $V^{\langle a \rangle}$ of $H^{1/4}$-normalized solutions of the mKdV equation \eqref{mkdv}, which are well controlled for 
an interval of time which increases without bound as $N\to\infty$.  These can be considered a weak analogue of the global solutions $u^{\langle aw \rangle}$ to NLS constructed in Section \ref{pseudo-sec}, but for mKdV,
and only for times $0 < t \ll \log N$ rather than $0 < t < \infty$. 
\begin{proof}[Conclusion of proof of Theorem \ref{defocus-mkdv}]
We now use scale invariance as in Section \ref{nls-sec} to construct $H^s$ solutions for $-1/4< s < 1/4$.

Let $0 < \delta \ll \eps \ll 1$ and $T > 0$ be arbitrary.  
As in Section \ref{nls-sec}, we shall find two solutions $u = \phi^{\langle a \rangle}, \phi^{\langle a' \rangle}$ to \eqref{mkdv} such that \eqref{0-small}, \eqref{0-close}, and \eqref{T-far} hold. 

Let $\lambda \gg 1$ be a large parameter to be chosen later.  Let $\phi^{\langle a \rangle}$ denote the function
$$ \phi^{\langle a \rangle}(t,x) := \lambda V^{\langle a \rangle}(\lambda^3 t, \lambda x).$$
Since $V^{\langle a \rangle}$ is a global smooth solution to mKdV, 
so also is $\phi^{\langle a \rangle}$. 

Similarly define
$$ \phi^{[a]}(t,x) := \lambda V^{[a]}(\lambda^3 t, \lambda x);$$
thus
$$\phi^{\langle a \rangle}(0,x) = \phi^{[a]}(0,x) = \lambda V^{[a]}(0,\lambda x) =
\lambda \sqrt{\frac{2}{3N}} \Re~  e^{iN\lambda x} u^{\langle aw \rangle}(0,\lambda x/\sqrt{3N}).$$
To estimate the $H^s_x$ norm of this function, we apply Lemma \ref{modulate} 
with $M=N\lambda$ and $\tau = N^{1/2}\lambda^{-1}$;
then $M\tau \equiv N^{3/2}\gg 1$ as $N\to\infty$.  Thus when $s \geq 0$
Lemma \ref{modulate} gives
$$ \| \phi^{\langle a \rangle}(0) \|_{H^s_x} \lesssim \lambda^{1/2+s} N^{s-1/4} \| u^{\langle aw \rangle}(0) \|_{H^1_x}.$$
For $-1/4 < s < 0$, we can still apply Lemma \ref{modulate} since $M \tau^{1 + s/K} \gg 1$ for sufficiently large $K$, the only difference being that the $H^1_x$ norm on the right must now be replaced by $H^K_x$.

If we thus define $\lambda$ by
\be{lam-def}
\lambda := N^{(1/4-s)/(1/2+s)}
\end{equation}
then by \eqref{uaw-bound} we have \eqref{0-small}:  
$\| \phi^{\langle a \rangle}(0) \|_{H^s_x}, \| \phi^{\langle a' \rangle}(0) \|_{H^s_x}\lesssim \eps$.
A similar argument using \eqref{soft-lip} instead of \eqref{uaw-bound} gives 
\eqref{0-close}:
$\| \phi^{\langle a \rangle}(0) - \phi^{\langle a' \rangle}(0) \|_{H^s_x} 
\lesssim \delta$.

Now we show \eqref{T-far}:  
$\sup_{0 \leq t < T} 
\| \phi^{\langle a \rangle}(t) - \phi^{\langle a' \rangle }(t) \|_{H^s_x} 
\gtrsim \eps$.
A routine scaling calculation shows
\be{scale}
\| \phi^{\langle a \rangle}(t) - \phi^{[a]}(t) \|_{H^s_x} \lesssim \lambda^{\max(s,0)+1/2} \| V^{\langle a \rangle}(\lambda^3 t) - V^{[a]}(\lambda^3 t) \|_{H^s_x}.
\end{equation}
Bounding the $H^s_x$ norm by the $H^{1/4}_x$ norm, we thus see from \eqref{09} that
\be{ee}
\| \phi^{\langle a \rangle}(t) - \phi^{[a]}(t) \|_{H^s_x} \lesssim \lambda^{\max(s,0)+1/2} \eps N^{-1+\eta}
\end{equation}
whenever $0 < t \ll \log N/\lambda^3$.  Applying \eqref{lam-def} and the hypothesis $s > -1/4$, we observe that the right-hand side is $\ll \eps$ if $\eta$ is chosen sufficiently small.

In particular, from Lemma \ref{modulate} we have
\be{phi-bound}
\| \phi^{\langle a \rangle}(t) \|_{H^s_x} \lesssim \eps.
\end{equation}

From \eqref{u-decohere} there exists a time $t_0>0$ depending on $a, a'$ (but not on $N$, $\lambda$) such that
$$ \| u^{\langle aw \rangle}(t_0) - u^{\langle a'w \rangle}(t_0) \|_{L^2_x} \gtrsim \eps.$$
Fix this $t_0$; we may choose $N$ so large that $t_0 \ll \log N$.  From \eqref{uaw-bound} and Lemma \ref{modulate} as before we thus have
$$ \|\phi^{[a]}(t_0/\lambda^3) - \phi^{[a']}(t_0/\lambda^3)\|_{H^s_x} \sim \lambda^{1/2+s} N^{s-1/4} \eps = \eps,$$
so by \eqref{ee} (and the remark immediately following) we have
$$ \|\phi^{\langle a \rangle}(t_0/\lambda^3) 
- \phi^{\langle a' \rangle}(t_0/\lambda^3)\|_{H^s_x} \sim \eps.$$
If we choose $N$ (and hence $\lambda$) large enough, we can make $t_0/\lambda^3 < T$, and so \eqref{T-far} follows.  This concludes the proof of Theorem \ref{defocus-mkdv}.
\end{proof}

\section{Energy estimates}\label{energy-sec}

In the last section we proved the ill-posedness of the mKdV equation below $H^{1/4}$.  However this proof is not as ``strong'' as the corresponding argument for NLS in Section \ref{nls}, because it did not rely on a family of global solutions to mKdV.   Instead, it relied on solutions which could only be controlled for very short times 
(roughly on the order of $\log N / N^c$ for some exponent $c>0$).  
Although this does suffice to disprove
uniform continuity of the solution map, it is not as satisfactory as the NLS argument, and does not give very good quantitative control on the ill-posedness.
 
In this section and the next, we rectify this shortcoming
by constructing a large family of global-in-time solutions to mKdV, 
similar to the global solutions $u = u^{\langle w \rangle}$ 
to NLS constructed in Section \ref{pseudo-sec}.  Unfortunately there seems to be no simple analogue of the pseudo-conformal transformation for mKdV, so our arguments will be more complicated.  On the other hand, these methods seem to be quite general and should be applicable to a wider class of nonlinear equations than those studied here.

In analogy with Section \ref{pseudo-sec} we expect our solutions $u$ to contain highly oscillating factors such as $e^{iNx}$.  These oscillations force
certain modifications in the estimates, but we use a {\it{lifting device}} 
to eliminate much of these difficulties.
This device is closely related to the method of ``slow and fast''
variables, which is frequently used in describing behavior in 
various asymptotic regimes.

This section and the next are devoted to the proof of:
\begin{theorem} \label{muchado}
Let $w\in C^\infty(\R)$ be real-valued and compactly supported.
Then for each sufficiently small $\eps>0$ and sufficiently large $\R\owns N<\infty$
there exist an exact solution $v^{\langle w\rangle}(x,t)$ and an approximate solution
$v^{[w]}(x,t)$ of the defocusing modified KdV equation $v_t+v_{xxx} = 6v^2v_x$,
defined for all $t\ge 2$, with the following properties:

The approximate solution takes the form
\begin{equation}
v^{[w]} = 2\eps N^{-1/2}t^{-1/2}\cos\big( \phi(x,t) \big)w(z)
\end{equation}
where
\begin{align*}
z &= N^{-1/2}t^{-1}\big(x + 3N^2t\big)
\\
\phi(x,t) &= -(-4x^3/27t)^{1/2}  
+ 6\eps^2 \log(t) N^{-1} (-x/3t)^{1/2}w^2(z),
\end{align*}
and satisfies
\begin{equation}
\|v^{[w]}(t)\|_{H^{1/4}}\sim\eps \hbox{ for all } t \geq 2. 
\end{equation}

The exact solution $v^{\langle w\rangle}$ is asymptotic to $v^{[w]}$
in the sense that
\begin{equation}
\|v^{\langle w\rangle}(t)-v^{[w]}(t)\|_{H^{1/4}}\lesssim \eps t^{-1/2}
\ \text{ for all } t\ge 2.
\end{equation}
\end{theorem}

Here the phase $\phi$ is approximately equal to $Nx+N^3t$, to leading order,
so heuristically 
\[
v^{\langle w\rangle}(t,x)\approx 2\eps N^{-1/2}t^{-1/2}\cos(Nx+N^3t)
w\big(N^{-1/2}t^{-1}(x + 3N^2t)\big).
\]

\subsection{Lifting }

Fix $N$ large.
We shall work in the cylinder $\R \times (\R/2\pi\Z) = \{ (y,\theta): y \in \R, \theta \in (\R/2\pi\Z) \}$.  We observe that we can embed the real line $\R$ into the cylinder by the map $x \mapsto (x/N^{1/2}, Nx)$; informally, this wraps $\R$ around the cylinder in a very tight spiral.

Our solution $u(t,x)$ shall be obtained by descent from a function $\tilde u(t,y,\theta)$ on the cylinder, via the transformation
\be{u-def}
u(t,x) := N^{-1/2} \tilde u(t, x/N^{1/2}, Nx).
\end{equation}
Indeed, observe from \eqref{u-def} that
$$ \partial_x u = N^{-1/2} (N^{-1/2} \partial_y + N \partial_\theta) \tilde u$$
and so if $\tilde u$ satisfies the PDE
\be{transformed-mkdv}
(\partial_t + (N^{-1/2} \partial_y + N \partial_\theta)^3) \tilde u
= 2 N^{-1} (N^{-1/2} \partial_y + N \partial_\theta) (\tilde u^3)
\end{equation}
on the cylinder, then $u$ will satisfy mKdV.

The lifting device expresses the highly oscillating function 
$e^{iNx}$ on the real line
as $e^{i\theta}$ on the cylinder, eliminating the dependence on $N$.
This will allow us to express certain nonstandard energy-type
estimates for functions of $x$ as more standard energy
estimates for functions of $(y,\theta)$. A disadvantage is that the dispersive term $(N^{-1/2} \partial_y + N \partial_\theta)^3$ has a large coefficient,
as $N\to\infty$.

%

We can also control $u$ in terms of $\tilde u$ by Sobolev norms,
by the following variant of the Sobolev trace lemma:

\begin{lemma}\label{sob-trace}  If $u$ and $\tilde u$ obey \eqref{u-def}, then
$$ \|u(t) \|_{H^{1/4}_x} \lesssim \| \tilde u(t) \|_{H^2_{y,\theta}}.$$
\end{lemma}

The argument below actually allows one
to lower $H^2$ to $H^{3/4+}$, which is consistent with the Sobolev trace 
lemma, but this yields no improvement in our application.  
If $\tilde u(t)$ has the special form $e^{ik\theta} a(y)$ for some small $k$ and smooth $a$ (e.g. $a \in H^2$), then one can also obtain corresponding lower bounds for the $H^{1/4}_x$ norm of $u(t)$ via Lemma \ref{modulate}.

\begin{proof}
It will suffice to prove the bound
\be{h1-bound}
\| u(t) \|_{L^2_x} \lesssim N^{-1/4} \| \tilde u(t) \|_{H^1_{y,\theta}}.
\end{equation}
Indeed, by applying the operator $\partial_x = N^{-1/2} \partial_y + N \partial_\theta$ to this estimate one obtains
$$ \| \partial_x u(t) \|_{L^2_x} \lesssim N^{-1/4} \| (N^{-1/2} \partial_y + N \partial_\theta) \tilde u(t) \|_{H^1_{y,\theta}} \lesssim N^{3/4} \| \tilde u \|_{H^2_{y,\theta}}$$
and then the claim follows by interpolation.  

It remains to prove \eqref{h1-bound}.  From Fubini's theorem we observe
$$
\| u(t) \|_{L^2_x}^2 = \int u(t,x)^2\ dx = C N^{-2} \int_0^{2\pi} 
\big| \sum_{k \in \Z} \tilde u(t, (k+\theta)/N^{3/2}, \theta)^2 \big| 
\ d\theta$$
where we have identified $\R/2\pi\Z$ with $[0,2\pi)$ in the obvious manner.  This it will suffice to prove the one-dimensional estimate
$$ N^{-3/2} \sum_{k \in \Z} |f(N^{-3/2}\theta)|^2 \lesssim \|f\|_{H^1_y}^2$$
for any function $f(y)$ and any $\theta \in [0,2\pi)$.

The left-hand side is bounded by
$$ \sum_{j \in \Z} \| f \|_{L^\infty([j,j+1])}^2,$$
which by the local Sobolev (or Poincar\'e) inequality is bounded by
$$ \sum_{j \in \Z} \| f \|_{H^1_y([j-1,j+2])}^2.$$
The claim follows.
\end{proof}

\subsection{Approximating the lifted evolution}

From the previous Lemma, we see that to construct global $H^{1/4}_x$ solutions $u$ to mKdV, 
it will suffice to construct global\footnote{
This $H^2_{y,\theta}$ norm may of course be reinterpreted
in the single $x$ coordinate, by replacing the Fourier multiplier
$1 + |\nabla_{y,\theta}|^2$ used to define $H^2_{y,\theta}$ by a multiplier 
$m_N(D)$, whose symbol $m_N(\xi)$ looks roughly like $m_N(kN + a) := N^{1/4} (1 + |k|)^2 (1 + N^{1/2} a)^2$ whenever $k$ is an odd integer and $|a| \leq N$.
Then the following analysis can be performed purely in the one-dimensional model using $m_N(D)$, and indeed this was our initial approach, but it is more complicated technically, mainly due to the need to develop a good Leibnitz rule for $m_N(D)$.
The lifting device allows us to work with more standard energy-type norms.
}
$H^2_{y,\theta}$ solutions $\tilde u$ to \eqref{transformed-mkdv}.  We will
do this in two steps,
constructing explicit global {\em approximate} solutions, then
modifying them to obtain exact solutions.
The following result, based on the energy method, 
asserts that any global approximate solution is asymptotic to
an exact solution, provided that the approximate solution
satisfies \eqref{transformed-mkdv} modulo a sufficiently small
remainder, as $t\to\infty$.

%

\begin{theorem}\label{perturb}  Let $0 < \eps \ll 1$ be a small number.  
Suppose that $\tilde v(t,y,\theta)=\tilde v_N(t,y,\theta)$ is a
one-parameter family of global real-valued smooth functions, 
rapidly decreasing in $y$, such that $\tilde v$ and the
error $\scripte=\scripte_N$ defined by
\be{transformed-mkdv-error}
\scripte:=(\partial_t + (N^{-1/2} \partial_y + N \partial_\theta)^3) \tilde v
- 2 N^{-1} (N^{-1/2} \partial_y + N \partial_\theta) (\tilde v^3) 
\end{equation}
obey the estimates
\begin{align}
\| \tilde v(t) \|_{C^3_{y,\theta}} &\lesssim \eps t^{-1/2}\label{v-infty}\\
\| \scripte(t) \|_{H^2_{y,\theta}} &\lesssim \eps t^{-\beta}   \label{e-x} 
\end{align}
for $t \ge 2$,
uniformly for all $N\ge 1$.
Suppose that 
\[
\beta >\frac32.
\]
Then, if $\eps$ is sufficiently small, 
there exists for each $N$ a global real-valued smooth solution 
$\tilde u (t, y , \theta )=\tilde u_N(t,y,\theta)$ 
to the transformed mKdV equation \eqref{transformed-mkdv} 
satisfying
\be{uv-diff}
\| \tilde u(t) - \tilde v(t) \|_{H^2_{y,\theta}} \lesssim \eps t^{1-\beta}
\end{equation}
for all $t \geq 2$, uniformly for $N\ge 1$.
In particular,
\[
t^{1/2}\| \tilde u(t) - \tilde v(t) \|_{H^2_{y,\theta}} \to 0
\text{ as } t\to\infty.
\]
\end{theorem}

The conditions \eqref{v-infty}, \eqref{e-x}
arise naturally in our construction of the approximate solution $\tilde v$ in 
the next section;  we will actually have 
$\|\scripte(t)\|_{H^2_{y,\theta}} \lesssim \eps t^{-2}(\log t)^C$.

\begin{proof}
We write \eqref{transformed-mkdv} as
$$ \tilde u_t = A\tilde u + D(\tilde u^3)$$
and \eqref{transformed-mkdv-error} as
$$ \tilde v_t = A\tilde v + D(\tilde v^3) + \scripte$$
where $A$ is the anti-self-adjoint constant-coefficient differential operator
$$ A := -(N^{-1/2} \partial_y + N \partial_\theta)^3$$
and $D$ is the constant coefficient 
vector field
$$ D := 2 N^{-3/2} \partial_y + 2 \partial_\theta,$$
whose coefficients are uniformly bounded for $N\ge 1$.

We pick a large time $T_0 \gg 1$ (which will eventually be set to infinity; all our bounds will be independent of $T_0$) and solve the backwards Cauchy problem
$$ \tilde u_t = A\tilde u + D(\tilde u^3); \quad \tilde u(T_0) = \tilde v(T_0)$$
on the region $2 \leq t \leq T_0$.  This has a global smooth solution (for the same reason that the Cauchy problem for mKdV has global smooth solutions; 
indeed one can foliate the cylinder into tightly wound copies of the real line on which the above equation is just a rescaled version of mKdV).  Writing $\tilde u = \tilde v + w$, we see that $w$ satisfies the difference equation
\be{w-eq}
w_t = Aw + D(w^3 + 3 \tilde v w^2 + 3\tilde v^2 w) - \scripte; \quad w(T_0) =0.
\end{equation}
We now introduce the energies
$$ \E_j(t) := \frac{1}{2} \int |\nabla_{y,\theta}^j w(t)|^2\ dx$$
for $j=0,1,2$.
Clearly $\E_j(T_0) = 0$ for $j=0,1,2$. We claim the estimates
\be{e-est}
\E_j(t) \leq C_0 \eps^2 t^{2-2\beta} \hbox{ for } j=0,1,2
\end{equation}
for all $2 \leq t \leq T_0$, where $C_0$ is a large absolute constant and assuming $\eps$ is sufficiently small (depending on $C_0$).  

To prove these estimates, we make the \emph{a priori} assumption that
\be{e-est-cont}
\E_j(t) \leq 2C_0 \eps^2 t^{2-2\beta} \hbox{ for } j=0,1,2
\end{equation}
for all $T \leq t \leq T_0$, and some $T \in [2, T_0]$.  We will then prove \eqref{e-est} for all $T \leq t \leq T_0$.  Since $w$ is smooth, this implies that the set of times $t$ for which \eqref{e-est} holds is both open and closed, and contains $T_0$.  From the continuity method we thus see that \eqref{e-est} will indeed hold for all $2 \leq t \leq T_0$.

We thus fix $T$ and assume \eqref{e-est-cont}.  In particular we have
\be{w-sobolev}
\| w(t) \|_{H^1_{y,\theta}} \lesssim \eps C_0^{1/2} t^{1-\beta}; \quad
\|w(t) \|_{C^1_{y,\theta}} \lesssim \| w(t) \|_{H^2_{y,\theta}} \lesssim \eps C_0^{1/2} t^{1-\beta}.
\end{equation}
By the hypothesis $\beta>3/2$, we have therefore the essential bound
\[ \|w(t) \|_{C^1_{y,\theta}}\lesssim\eps t^{-1/2}.\]

We differentiate $\E_0$ to obtain 
$\E'_0(t) = \int ww_t$, and substitute \eqref{w-eq} for $w_t$.  Since $A$ is anti-self-adjoint and commutes with $\nabla_{y,\theta}$, its contribution to $\E_0'(t)$ vanishes, and we obtain
$$ \E_0'(t) = \int w \cdot
 \big(D(w^3 + 3 \tilde v w^2 + 3\tilde v^2 w) - \scripte\big).$$
We expand out the cubic terms, using the Leibnitz rule.  Any term of the form $\int w Dw \tilde v \tilde v$ 
can be rewritten as $-\int w w \tilde v D\tilde v$ using the identity $w Dw = \frac{1}{2}D(w^2)$ and integration by parts.  Thus we obtain a finite sum
of integrals, in each of which there are least two factors of $w$ 
on which no derivatives fall.
We then use Cauchy-Schwarz to obtain
$$ |\E_0'(t)| \lesssim \| w\|_{L^2_{y,\theta}}^2 
\big(\|w\|_{C^1_{y,\theta}} + \|\tilde v\|_{C^1_{y,\theta}}\big)^2 
+ \|w\|_{L^2_{y,\theta}} \| \scripte \|_{L^2_{y,\theta}}.$$
Applying \eqref{w-sobolev}, \eqref{v-infty}, \eqref{e-x} we obtain
$$ |\E_0'(t)| \lesssim C_0^2 \eps^4 
t^{2-2\beta}(t^{-1/2})^2  
+ C_0^{1/2} \eps^2 t^{1-\beta}t^{-\beta}.$$
Integrating this (using $\E_0(T_0) = 0$) we thus obtain 
$\E_0(t) \leq C_0 \eps^2 t^{2-2\beta}$ 
as desired, provided $C_0$ is sufficiently large and $\eps$ sufficiently small (depending on $C_0$).

The higher-order quantities $\E_1,\E_2$ are handled in the same way.
Consider first $\E_1$.
By arguing as before (and noting that $A$ commutes with $\nabla_{y,t}$) we have
$$ \E_1'(t) = \int \nabla_{y,t} w \cdot \nabla_{y,t} (D(w^3 + 3 \tilde v w^2 + 3\tilde v^2 w) - \scripte).$$
The most dangerous terms are those involving a $\nabla_{y,t}w \cdot D \nabla_{y,t}w$ factor, but by using the identity
$$ \nabla_{y,t} w \cdot D \nabla_{y,t} w = \frac{1}{2} D(|\nabla_{y,t} w|^2)$$
and integrating by parts we may transfer $D$ to another factor,
as we did in analyzing $E_0$.  Thus no factor $w$ with two derivatives on it 
will remain. 
Unfortunately, there may remain a factor $\tilde v$ 
with two derivatives on it, 
but if that is the case then all but one of the factors of $w$ 
will have no derivatives. Applying H\"older, we then obtain
$$ |\E'_1(t)| \lesssim \| w \|_{H^1_{y,t}}^2 
(\big \|w\|_{C^1_{y,t}} +  \|\tilde v\|_{C^1_{y,t}} \big)^2
+ \|w\|_{H^1_{y,t}} \|w\|_{L^2_{y,t}}
\|w\|_{L^\infty_{y,t}}  \|\tilde v\|_{C^2_{y,t}} 
+ \| w\|_{H^1_{y,t}} \| \scripte\|_{H^1_{y,t}}.$$

To majorize $\E_2$, we differentiate with respect to $t$ and argue as before,
obtaining
$$ \E'_2(t) = \int \nabla_{y,t}^2 w \cdot (D \nabla_{y,t}^2 (w^3 + 3 \tilde v w^2 + 3\tilde v^2 w) - \scripte).$$
We again apply the Leibnitz rule.  Again the most dangerous terms are those 
with $\nabla_{y,t}^2 w \cdot D \nabla_{y,t}^2 w = \frac{1}{2} D(|\nabla_{y,t}^2 w|^2)$, but each such term may be rewritten by distributing
 the $D$ to a factor of $w$ on which no other derivatives fall.
Thus no factor of $w$ will carry more than two derivatives.
Since there are at most five derivatives in any of the terms, for any term involving two factors of $\nabla_{y,t}^2 w$, all other terms carry
derivatives of at most first order.  
From this and Cauchy-Schwarz we see that
$$ |\E'_2(t)| \lesssim \| w\|_{H^2_{y,t}}^2 (\| w\|_{C^1_{y,t}} 
+ \| \tilde v \|_{C^1_{y,t}})^2 + \| w \|_{H^2_{y,t}} \|w\|_{H^1_{y,t}} 
(\| w\|_{C^1_{y,t}} + \| \tilde v \|_{C^3_{y,t}})^2 + \| w \|_{H^2_{y,t}} \| \scripte\|_{H^2_{y,t}}.$$

From these two differential inequalities we deduce as for $\E_0$ that 
$\E_j(t) \leq C_0 \eps^2 t^{2-2\beta}$ for $j=1,2$.
This concludes the proof of \eqref{e-est}.

Finally, we need to remove the restriction $t < T_0$.  This can be achieved by letting $T_0 \to +\infty$, and taking a weak limit in $H^2_x$ of the functions $w(t) = w^{(T_0)}(t)$ (which is thus strongly convergent in $H^{2-\eps}_x$, by Rellich embedding); observe that all the above bounds were independent of $T_0$.  One then obtains in the limit a new function $w^{(+\infty)}(t)$ which obeys \eqref{w-eq} and \eqref{e-est} (and hence \eqref{uv-diff}, if 
$\tilde u := \tilde v + w^{(+\infty)}$) for all $t \geq 2$.  
We omit the details.
\end{proof}

\section{Construction of the approximate solution}\label{approx-sec}

In this section we
construct solutions $\tilde v$ to \eqref{transformed-mkdv-error} which obey the bounds \eqref{v-infty}, \eqref{e-x}. 
For this task it is more convenient 
to work with the original equation, rather than in the $(y,\theta)$ 
variables, mainly because $\partial_x^3$ is often more convenient
to work with than 
$(N^{-1/2} \partial_y + N \partial_\theta)^3$.  
We will begin by constructing a family of approximate real solutions 
$v = v^{[\varphi_1]}$ to mKdV. These depend on a bump function 
$\varphi_1$ which is smooth and rapidly decaying, but is otherwise
arbitrary.
We introduce the Ansatz
\begin{equation} \label{vdefinition}
v := v_{-3} + v_{-1} + v_1 + v_3
\end{equation}
where the $v_k$ are complex functions depending upon $\varphi_1$ 
(which one should think of as 
oscillating essentially like $e^{iNkx}$ in space) obeying $v_{-k} = \overline{v_k}$.  The dominant terms will be $v_{\pm 1}$; for instance, their $L^\infty$ norm will be $O(N^{-1/2} t^{-1/2})$,
whereas $v_{\pm 3}$ will be $O(N^{-7/2} t^{-3/2})$ in $L^\infty$.  
However, the correction terms $v_{\pm 3}$ will be necessary in order to obtain the $t^{-2} \log^C(t)$ decay \eqref{e-x} for the error term $\scripte$;
without these terms, the error turns out to decay only like 
$t^{-1} \log^C(t)$.
We shall assume $t \geq 2$ throughout to avoid the (artificial) singularity at $t=0$.

We now construct an approximate
solution $v = v_{-3} + v_{-1} + v_1 + v_3$ to mKdV which has the properties stated in the previous section. In fact we shall construct a family of such solutions $v$ which depend on an arbitrarily chosen bump function $\varphi_1(z)$.  

With $v$ of the general form \eqref{vdefinition},
$$ v_t + v_{xxx} = 2(v^3)_x + E$$
where the error $E = \sum_{k \hbox{ odd}:\,\, |k| \leq 9} E_k$ is given by
$$ E_k := (\partial_t + \partial_{xxx}) v_k - 2 \partial_x \sum_{k_1 + k_2 + k_3 = k} v_{k_1} v_{k_2} v_{k_3}.$$
In particular we have $E_{-k} = \overline{E_k}$.  We note in particular that
\bas
E_1 &= (\partial_t + \partial_{xxx}) v_1 - 6 \partial_x(|v_1|^2 v_1) + \tilde E_1 \\
E_3 &= (\partial_t + \partial_{xxx}) v_3 - 2 \partial_x(v_1^3) + \tilde E_3 \\
E_k &= \tilde E_k\ \text{for } k=5,7,9
\end{align*}
where the $\tilde E_k$ are linear combinations of expressions of the form
$\partial_x(v_{k_1} v_{k_2} v_{k_3})$, where $k_1+k_2+k_3=k$ and at least one of the $k_1$, $k_2$, $k_3$ is equal to $\pm 3$.  Heuristically,
the $\tilde E_k$ terms will be negligible because $v_3$ is much smaller than $v_1$; 
most of the work will arise in controlling the dominant terms in $E_1$, 
and to a lesser extent in $E_3$. 

We introduce the corresponding functions of the $(y,\theta)$ variables
\bas
\tilde v(t,y,\theta) &:= \sum_{k=-3,-1,1,3} N^{1/2} e^{ik\theta} e^{-iN^{3/2}ky} v_k(t, N^{1/2} y)\\
\scripte(t,y,\theta) &:= \sum_{k=-9,\ldots,9} N^{1/2} e^{ik\theta} e^{-iN^{3/2}ky} E_k(t, N^{1/2} y),
\end{align*}
and observe (from the fact that 
$\partial_x = N^{-1/2} \partial_y + N \partial_\theta$ 
annihilates $e^{ik\theta} e^{-iN^{3/2}ky}$) 
that $\tilde v$ and $\tilde \scripte$ 
obey \eqref{transformed-mkdv-error}.  
Thus we will be able to invoke Theorem \ref{perturb} provided that we are able
to construct $v_k$ (and hence $E_k$) obeying the estimates
\begin{align}
\| e^{-iN^{3/2} ky} v_k(t, N^{1/2} y) \|_{C^3_y} &\lesssim \eps N^{-1/2} t^{-1/2} \label{v-trans}\\
\| e^{-iN^{3/2} ky} E_k(t, N^{1/2} y) \|_{H^2_y} &\lesssim \eps N^{-1/2} t^{-2} \log^C t \label{E-trans}
\end{align}
for all $t \geq 2$ and all $k$ for which the above make sense.  

Introduce the coordinate
\begin{equation}  \label{zdefn}
z := N^{-1/2}t^{-1}x + 3N^{3/2} = N^{-1/2}t^{-1}\big(x + 3N^2t \big).
\end{equation}
We will work partly in coordinates $(x,t)$, and partly in
coordinates $(z,t)$;
we will always work in a region in which $z$ is uniformly bounded,
so that 
\begin{equation}
\frac{-x}{3t} = N^2(1-\tfrac13 N^{-3/2}z) = N^2+O(N^{1/2}).
\end{equation}
Thus we have $x = - 3N^2 t + O(N^{1/2} t)$, and in particular $x$ 
is always negative\footnote{This reflects the fact that solutions to the 
Airy equation tend to propagate rapidly to the left, especially if the 
solution is high frequency as is the case here.} and $|x| \sim N^2 t$.   
Thus fractional powers of $-x$ are well-defined.

We set
\be{vk-def}
v_k = \eps^k N^{(2-3|k|)/2} t^{-|k|/2} e^{ik\phi(t,x)} \varphi_k(z),
\end{equation}
 for various real-valued bump functions $\varphi_k$ to be chosen later, and a 
real phase function $\phi(t,x)$ which is chosen to satisfy a naturally
arising eikonal-type equation.  
Of course, we choose $\varphi_{-k} = \varphi_k$ in order to have 
$v_k = \overline{v_{-k}}$.

Suppose $v_k$ takes the form \eqref{vk-def}. From
\bas
\partial_t e^{i k\phi} &= i k\phi_t e^{ik\phi} \\
\partial_x e^{i k\phi} &= i k\phi_x e^{ik\phi} \\
\partial_{xx} e^{i k\phi} &= (-k^2\phi_x^2 + i k\phi_{xx}) e^{ik\phi} \\
\partial_{xxx} e^{i k\phi} &= (-ik^3\phi_x^3 - 3k^2\phi_x \phi_{xx} + i k\phi_{xxx}) e^{ik\phi}
\end{align*}
and
\bas
\partial_t \varphi(z) &= -N^{-1/2} x t^{-2} \varphi '(z) \\
\partial_x \varphi(z) &= N^{-1/2} t^{-1} \varphi'(z) \\
\partial_{xx} \varphi(z) &= N^{-1} t^{-2} \varphi''(z) \\
\partial_{xxx} \varphi (z) &= N^{-3/2} t^{-3} \varphi'''(z) 
\end{align*}
follows the fundamental formula
\be{ident}
\begin{split}
(\partial_t + \partial_{xxx}) v_k
=
\eps^k N^{1-3|k|/2} t^{-|k|/2}
e^{ik\phi}
\cdot
\Big[
&(i k\phi_t - ik^3 \phi_x^3 - 3 k^2 \phi_x \phi_{xx} - \frac{|k|}{2} t^{-1} + i k \phi_{xxx} )  \varphi_k(z)\\
&+ (-\frac{x}{t} - 3k^2 \phi_x^2 + 3i k\phi_{xx})N^{-1/2}t^{-1}\varphi'_k(z) \\
&+ 3i k \phi_x  N^{-1} t^{-2}  \varphi''_k(z) \\
&+ N^{-3/2}t^{-3} \varphi'''_k(z)
\Big].
\end{split}
\end{equation}
Terms that decay like $t^{-5/2}$ or better will turn out to be negligible.  
When $k=1$, this will include the terms involving $\varphi''_1$, $\varphi'''_1$, $\phi_{xx} \varphi'_1$, or $\phi_{xxx}$; when $k=3$, this will include 
all terms except the very first two, $ik\phi_t-ik^3 \phi_x^3$.

We have not yet specified what the phase function $\phi$ is.  Before proceeding with the detailed estimation of the terms in \eqref{ident}, we briefly indicate the heuristic considerations which lead naturally to the choice of this function. 
As initial approximations to solutions of mKdV,
we begin with solutions
$\int e^{ix\xi+it\xi^3}h(\xi)\,d\xi$ 
of the Airy equation $(\partial_t + \partial_{xxx})u=0$.
Assuming that $h\in C_0^\infty$ is supported in a compact
subset of $(0,\infty)$, the stationary phase method
gives the leading-order asymptotics as $t\to+\infty$ to be
$ce^{i\Phi(x,t)}t^{-1/2}\tilde h((-x/3t)^{1/2})$
where 
\be{Phidef}
\Phi(x,t) = -(-4x^3/27t)^{1/2} 
\end{equation}
and $\tilde h$ is another bump function, $\tilde h(\xi)= h(\xi)\xi^{-1/2}$.
We therefore take \eqref{Phidef} as an initial approximation to $\phi$.

Taking $\phi=\Phi$ in \eqref{ident}, one finds that the right-hand
side is $O(t^{-5/2})$, as is desired. However,
$(\partial_t+\partial_{xxx})v_1-6\partial_x(|v_1|^2v_1)$
is larger. Indeed, the main term of $6\partial_x(|v_1|^2v_1)$
is $6i\eps^3 N^{-3/2}t^{-3/2}\Phi_x\varphi_1^3e^{i\phi}$; see \eqref{mainNL} below. 
Thus we set $\phi = \Phi+\psi$ and solve for $\psi$
by setting the main new term resulting from the replacement of $\phi=\Phi$
by $\phi=\Phi+\psi$ in \eqref{ident} equal to 
this main term of $6\partial_x(|v_1|^2v_1)$.
This gives 
\[
i(\partial_t-3\Phi_x^2\partial_x)\psi
\cdot \big(\eps N^{-1}t^{-1/2} \big)e^{i\phi}
=
6i\eps^3 N^{-3/2}t^{-3/2}\Phi_x\varphi_1^3e^{i\phi},
\]
whence
\begin{equation}
(\partial_t-3\Phi_x^2\partial_x)\psi
=
6\eps^2t^{-1}N^{-1}\Phi_x\varphi_1^2.
\end{equation}
When rewritten in the coordinates $(z,t)$, $\partial_t-3\Phi_x^2\partial_x$ becomes just $\partial_t$. Since $\Phi_x=(-x/3t)^{1/2}$ is a function of $z$ alone,
the latter equation may be solved explicitly: 
\[
\psi= 6\eps^2 \log(t)N^{-1}\Phi_x\varphi_1^2.
\]
This phase correction should be compared with \eqref{explicit}; it is relatively small compared to the dominant term $\Phi$ of the phase 
(which is $O(N^3t)$ compared to $O(\log t)$ for the phase correction).

We therefore define
\be{phi-def}
\phi (t,x) := \Phi(t,x) + \eps^2 \tilde \phi(z) \log(t)
\end{equation}
where
\begin{equation}
\tilde \phi(z) := 6N^{-1} \Phi_x \varphi^2_1(z) = 
6 (\frac{-x}{3N^2t})^{1/2} \varphi^2_1(z) 
= 6(1 - \tfrac13N^{-3/2}z)^{1/2} \varphi^2_1(z).
\end{equation}
$\varphi_1$ will be an arbitrary smooth function with compact support.
However, it remains to specify $\varphi_3$, which must be chosen to satisfy
an equation (see \eqref{varphi3-def} below)
in order that $E_3$ will be sufficiently small for our purpose.


A convention will simplify the notation.
We write $f(t,x)=\scripto(N^\alpha t^\beta \log(t)^\gamma)$
to mean that uniformly for all $t\ge 2$ and $(x,t)$ in the support of $\phi_1(z)$,
$$ |\partial_t^a \partial_x^b f(t,x)| 
\leq C_{f,a,b} N^\alpha t^\beta \log(t)^\gamma t^{-a} (N^{1/2} t)^{-b}$$
for all $a,b\ge 0$.
When $\alpha=\beta=\gamma=0$ we will often write $f=\scripto(1)$.  Observe in particular that $\varphi(z) = \scripto(1)$ for any $C^\infty$ function $\varphi$.

By Taylor expansion, 
$$ (\frac{-x}{3N^2 t})^\alpha = (1 - \tfrac13 N^{-3/2} z)^\alpha 
= 1 + \scripto(N^{-3/2}) = \scripto(1)$$
for any $\alpha \in \R$.  In particular we have
\bas
\Phi_t &= (\frac{-x^3}{27 t^3})^{1/2} = N^3 + \scripto(N^{3/2}) \\
\Phi_x &= (\frac{-x}{3 t})^{1/2} = N + \scripto(N^{-1/2}) \\
\Phi_{xx} &= -(\frac{-1}{12xt})^{1/2} 
= -\frac{1}{6} N^{-1} t^{-1} + \scripto(N^{-5/2} t^{-1}) \\
\Phi_{xxx} &= (\frac{-1}{48 x^3 t})^{1/2} = \scripto(N^{-3} t^{-2}). 
\end{align*}
If we then add in the phase correction $\eps^2 \tilde \phi(z) \log(t) 
= \scripto(\log(t))$ we obtain the following estimates for $\phi$:
\begin{align}
\phi_t &= (\frac{-x^3}{27 t^3})^{1/2} + \eps^2 \tilde \phi(z) t^{-1} - \eps^2 N^{-1/2}  x t^{-2} \tilde \phi'(z) \log(t) 
= N^3 + \scripto(N^{3/2}) \label{phit}\\
\phi_x &= (\frac{-x}{3 t})^{1/2} + \eps^2 N^{-1/2} t^{-1} \log(t) \tilde \phi'(z) = N + \scripto(N^{-1/2}) \label{phix}\\
\phi_{xx} &= -(\frac{-1}{12xt})^{1/2} + \eps^2\scripto(N^{-1} t^{-2} \log(t)) 
= \scripto(N^{-1} t^{-1}) \label{phixx}\\
\phi_{xxx} &= \scripto(N^{-3/2} t^{-2}) .\label{phixxx}
\end{align}

The functions $v_k$ and $E_k$ will be linear combinations of expressions of the form $e^{ik\phi} \scripto(1)$.  
To prove the desired estimates \eqref{v-trans}, 
\eqref{E-trans} we use the following lemma.

\begin{lemma}\label{xk}  Let $k = O(1)$ be an integer and $t \geq 2$, and suppose that $f_k(t,x) = e^{ik\phi} \scripto(1)$ 
is supported where $|z|=|N^{-1/2}t^{-1}x+3N^{3/2}|$ is bounded by some fixed constant.
Let $g(y)$ denote the function 
$$ g(y) := e^{-iN^{3/2} ky} f_k(t, N^{1/2} y).$$
Then
\bas
\| g(y) \|_{C^3_y} &\lesssim 1\\
\| g(y) \|_{H^2_y} &\lesssim t^{1/2}. 
\end{align*}
\end{lemma}

\begin{proof}
The $H^2$ norm bound follows directly from the $C^3$ bound and the support
hypothesis on $f$.
Since for any $n$, $\partial_y^n\scripto(1)= O(t^{-n})=O(1)$,
it suffices to verify that
$$ \partial_y^j (\phi(t,N^{1/2}y) - N^{3/2} y) = O(1)$$
on the support of $\varphi_1(y)$ for all $j=1,2,3$.
But this follows directly from \eqref{phix}, \eqref{phixx}, \eqref{phixxx}.
\end{proof}

From the above lemma we see immediately that $v$ obeys \eqref{v-trans}. 
To prove \eqref{E-trans}, it will suffice to show that 
$E_k = \eps \scripto(N^{-1/2} t^{-5/2} \log^C(t)) e^{ik\phi} $ 
for all $k$.  (Of course, any term in $E_k$ which decays even faster in $t$ or has more powers of $\eps$ and $N^{-1}$ will also be acceptable).  We may of course restrict our attention to positive $k$ since $E_{-k} = \overline{E_k}$.

We begin by computing certain expressions involving the phase $\phi$
which appear in \eqref{ident}.
\begin{lemma}\label{phi2-est} We have
\begin{align}
i \phi_t - i \phi_x^3 &= 6iN^{-1} t^{-1}\phi_x \eps^2 \varphi_1^2(z) 
+ \scripto(t^{-2} \log^2(t)) \label{a-phi1}\\
3i \phi_t - 27 i \phi_x^3 &= -24 i (\frac{-x}{3 t})^{3/2} 
+ \scripto(N^{3/2} t^{-1} \log (t)) \nonumber \\
&= -24 i N^3 (1 - \tfrac13 N^{-3/2}z)^{3/2}
+ \scripto(N^{3/2} t^{-1} \log (t)) \label{a-phi4}\\
- 3 \phi_x \phi_{xx} - \frac{1}{2t} &= \scripto(t^{-2} \log (t)) 
\label{a-phi2}\\
-\frac{x}{t} - 3 \phi_x^2 &= \scripto(N^{1/2} t^{-1} \log (t)) \label{a-phi3}
\end{align}
\end{lemma}

\begin{proof}
We begin with \eqref{a-phi1}, \eqref{a-phi4}.  From \eqref{phix} we have
$$ \phi_x^3 = (\frac{-x}{3 t})^{3/2}
+ 3 (\frac{-x}{3 t}) \eps^2 N^{-1/2} t^{-1} \log(t) \tilde \phi'(z)
+ \scripto(t^{-2} \log^2(t)) $$
while from \eqref{phit} we have
$$ \phi_t = (\frac{-x}{3 t})^{3/2} + \eps^2 \tilde \phi(z) t^{-1} - \eps^2 N^{-1/2} x t^{-2}  \log(t) \tilde \phi'(z).$$
Meanwhile, we have
\bas
6iN^{-1}t^{-1} \phi_x \eps^2 \varphi_1^2(z) 
&= 6iN^{-1}t^{-1}  \Phi_x \eps^2 \varphi_1^2(z) 
+ \scripto(N^{-1/2} t^{-2} \log(t)) \\
&= i \eps^2 t^{-1}\tilde \phi(z) 
+ \scripto(N^{-1/2} t^{-2} \log(t)) 
\end{align*}
The claims \eqref{a-phi1}, \eqref{a-phi4} follow.

Now we prove \eqref{a-phi2}.  From \eqref{phix} and \eqref{phixx} we have
$$ \phi_x = (\frac{-x}{3 t})^{1/2} + \scripto(N^{-1/2} t^{-1} \log(t)) ; 
\quad \phi_{xx} = -(\frac{-1}{12xt})^{1/2} + \scripto(N^{-1} t^{-2} \log(t)) $$
and hence
$$ \phi_x \phi_{xx} = -\frac{1}{6t} + \scripto(t^{-2} \log(t)) $$
which is \eqref{a-phi2}.  Likewise
$$ \phi_x^2 = \frac{-x}{3 t} + \scripto(N^{1/2} t^{-1} \log(t)) $$
which is \eqref{a-phi3}.
\end{proof}

Using this Lemma, \eqref{ident}, \eqref{phixx}, and \eqref{phixxx} 
we can now expand 
\bas
(\partial_t + \partial_{xxx}) v_1
=
\big({6}i{N^{-1}} t^{-1}\phi_x 
\eps^2 \varphi_1^2(z) 
&+ \scripto(t^{-2} \log(t))\big) t^{-1/2} 
\eps N^{-1/2} e^{i\phi} \varphi_1(z)\\
&+ \scripto(N^{1/2} t^{-1} \log(t))  
\eps N^{-1} t^{-3/2} e^{i\phi} \varphi_1'(z) \\
&+ \scripto(t^{-5/2}) \eps N^{-1/2} e^{i\phi} \varphi_1''(z) \\
&+ t^{-7/2} \eps N^{2} e^{i\phi} \varphi_1'''(z),
\end{align*}
which simplifies to
$$
(\partial_t + \partial_{xxx}) v_1
=
6i \eps^3 N^{-3/2} t^{-3/2}  \phi_x e^{i\phi} \varphi_1^3(z) 
+ \eps \scripto(N^{-1/2} t^{-5/2} \log(t)) e^{i\phi} .
$$
Meanwhile, we have
\begin{equation} \label{mainNL}
6\partial_x(|v_1|^2 v_1) 
= 6 \eps^3 N^{-3/2} t^{-3/2} i \phi_x e^{i\phi} \varphi_1^3(z) 
+ 18 \eps^3 N^{-2} t^{-5/2} e^{i\phi} \varphi_1^2(z) \varphi'_1(z).
\end{equation}
As foreshadowed in our heuristic derivation of $\phi$, the two leading terms here match.  From the definition of $E_1$ we thus have
$$ E_1 = \tilde E_1 + \eps \scripto(N^{-1/2} t^{-5/2} \log(t)) e^{i\phi} .$$
The last term is of the desired form; we will see below that 
$\tilde E_1$ is also.

Before verifying this, we turn to $E_3$.  
We begin by expanding $(\partial_t + \partial_{xxx}) v_3$.  
Using \eqref{ident}, estimating the main terms using \eqref{a-phi4}, 
and majorizing all the other terms crudely by 
\eqref{phix}, \eqref{phixx}, and \eqref{phixxx}, 
we can write this expression as 
\bas
(\partial_t+\partial_x^3)v_3 
=
\Big(-24 i (\frac{-x}{3 t})^{3/2} 
&+ \scripto(N^{3/2}t^{-1} \log(t)) \Big) 
\eps^3 N^{-7/2} t^{-3/2} e^{3i\phi} \varphi_3(z)\\
&+ (N^2 \scripto(1)) \eps^3 N^{-4} t^{-5/2} e^{3i\phi} \varphi_3'(z) \\
&+ \scripto(1)\eps^3 N^{-7/2} t^{-7/2} e^{3i\phi} \varphi_3''(z) \\
&+ \eps^3 N^{-5} t^{-9/2}  e^{3i\phi} \varphi_3'''(z),
\end{align*}
which simplifies to
$$
(\partial_t+\partial_x^3)v_3 
=
-24 i (\frac{-x}{3 t})^{3/2} 
\eps^3 N^{-7/2} t^{-3/2} e^{3i\phi} \varphi_3(z)
+ \eps \scripto(N^{-1/2} t^{-5/2} \log(t)) e^{3i\phi} .$$

Meanwhile, we have
\begin{align*} 
2\partial_x( v_1 ^3 ) 
&= 6 \eps^3 N^{-3/2} t^{-3/2} i \phi_x e^{3i\phi} \varphi_1^3(z) 
+ 6 \eps^3 N^{-2} t^{-5/2} e^{3i\phi} \varphi_1^2(z) \varphi'_1(z)
\\
&=
6i\eps^3N^{-3/2}t^{-3/2}\varphi_1^3
\big[
(\frac{-x}{3t})^{1/2} + \scripto(N^{-1/2}t^{-1}\log(t))
\big] e^{3i\phi}
+ \eps^3 \scripto(N^{-2} t^{-5/2}) e^{3i\phi} 
\\
&=
6i\eps^3N^{-3/2}t^{-3/2}\varphi_1^3
(\frac{-x}{3t})^{1/2}
e^{3i\phi}
+ \eps^3 \scripto(N^{-2} t^{-5/2}\log(t)) e^{3i\phi} .
\end{align*}
If we equate the first term of this last line with
the leading term of $(\partial_t+\partial_x^3)v_3$, then we obtain
the relation
\begin{equation}
-24i (\frac{-x}{3t})^{3/2}
\eps^3 N^{-7/2}t^{-3/2}\varphi_3
=
6i\eps^3N^{-3/2}t^{-3/2}\varphi_1^3
(\frac{-x}{3t})^{1/2}
(1+ \scripto(N^{-3/2})).
\end{equation}
Therefore if we define
\begin{equation} \label{varphi3-def}
\varphi_3(z) = -\tfrac14 \varphi_1^3(z) N^2 (\frac{-x}{3t})^{-1}
= - \frac{1}{4} \varphi_1^3(z) (1 - \frac{1}{3} N^{-3/2} z)^{-1} = \scripto(1),
\end{equation}
then 
\[
(\partial_t+\partial_x^3)v_3
-2\partial_x(v_1^3)
= \eps^3\scripto(N^{-1/2}t^{-5/2}\log(t))e^{3i\phi}.
\]
In particular, the support of $\varphi_3$ is a subset of
the support of $\varphi_1$, hence is bounded in the $z$ coordinate,
uniformly in $t$. Consequently all $v_j$ and $E_j$ share this
same uniform support property.

From the definition of $E_3$ we thus have
$$ E_3 = \tilde E_3 + \eps \scripto(N^{-1/2} t^{-5/2} \log(t)) e^{3i\phi} .$$
In light of these estimates, it will thus suffice to control the minor errors 
$\tilde E_k$, i.e.\ to show that
$$ \tilde E_k = \eps \scripto(N^{-1/2} t^{-5/2} \log(t)) e^{3ik\phi} $$
for $k = 1,3,5,7,9$.  Expanding out $\tilde E_k$, it thus suffices to show that
$$ (v_{k_1})_x v_{k_2} v_{k_3} 
= \eps^3 \scripto(N^{-1/2} t^{-5/2} \log(t)) e^{3i(k_1+k_2+k_3)\phi} $$
for all $k_1, k_2, k_3 \in \{-3,-1,1,3\}$ with at least one of 
$k_1, k_2, k_3$ equal to $\pm 3$.  But this follows from the estimates
\bas
v_{\pm 1} &= \eps \scripto(N^{-1/2} t^{-1/2}) e^{\pm i\phi} \\
v_{\pm 3} &= \eps \scripto(N^{-7/2} t^{-3/2}) e^{\pm 3i\phi} \\
(v_{\pm 1})_x &= \eps \scripto(N^{1/2} t^{-1/2}) e^{\pm i\phi} \\
(v_{\pm 3})_x &= \eps \scripto(N^{-5/2} t^{-3/2}) e^{\pm 3i\phi} ,
\end{align*}
which come from \eqref{vk-def} and \eqref{phix}.  
(Indeed, there is substantial room to spare, in terms of powers of $N$.)  

This completes the proof of \eqref{E-trans}.  Thus all the conditions of Theorem \ref{perturb} are obeyed for this choice $v = v^{[\varphi_1]}$ of approximate solution.  Applying this theorem followed by Lemma \ref{sob-trace}, we see that we can construct global $H^{1/4}$ solutions $u = u^{\langle \varphi_1\rangle}$ depending on an initial choice of bump function $\varphi_1$ which can be approximated in $H^{1/4}$ as $t \to +\infty$ by an explicit function $v$ given by the above Ansatz.  These solutions $u^{\langle\varphi_1\rangle}$ are closely analogous to the global solutions $u^{\langle w \rangle}$ to NLS constructed in Section \ref{pseudo-sec}.  As with NLS, the logarithmic factor in the phase leads to a proof of Theorem \ref{defocus-mkdv} which is closely analogous to the proof of Theorem \ref{defocus-nls} (and yields similar quantitative control of the nature of the ill-posedness).  Since ill-posedness has already been established by an alternative argument, we omit the details. 

However, one small comment is needed.  In establishing ill-posedness in 
$H^s$ for $s<0$,
a problem arises close to the scaling threshold $s=-1/2$ 
if the solutions $u^{\langle \varphi_1 \rangle}$ have a substantial low frequency component, as this will not scale favorably.  This can be ruled out
by observing that the function $\tilde v(t,y,\theta)$ constructed earlier in this section has the symmetry
$$ \tilde v(t,y,\theta + \pi) = - \tilde v (t,y,\theta)$$
(because all the integers $k$ in the summation are odd).  An inspection of the proof of Theorem \ref{perturb} reveals that $\tilde u$ must also have this symmetry (since it is preserved by the flow \eqref{transformed-mkdv}).  
Thus when 
$\tilde u$ is expanded in Fourier series in the angular variable $\theta$,
Fourier components $e^{ik\theta}$ with nonzero coefficients arise only for
odd $k$.  In particular, there is no zero Fourier mode.  From this and a variant of Lemma \ref{modulate} one can show that the solutions $u^{\langle \varphi_1 \rangle}$ will be extremely small at the frequency origin (especially if the 
$H^2_{y,\theta}$ control on $\tilde u$ is improved to $H^l$ for sufficiently
large $l$, in order to control higher Sobolev norms).  
Again, we omit the details. 

\section{The Miura transform}\label{miura-sec}

In this section we review the Miura transform relating solutions of
defocusing mKdV \eqref{mkdv} to real solutions of KdV \eqref{kdv}, and
show how this transform, combined with Theorem \ref{defocus-mkdv},
gives Theorem \ref{kdv-illp}.  In the next section we will introduce a generalization of the Miura transform (related to the Gardner transform) which will be used to prove Theorem \ref{kdv-endpoint}.

The \emph{Miura transform} $M$ is defined by
$$ M(v) := v_x + v^2.$$
Observe that if $v$ is a smooth real-valued solution to the mKdV equation \eqref{mkdv}, then $u = M(v)$ is a smooth real-valued solution to KdV equation \eqref{kdv}.  Indeed:
\begin{align*}
(\partial_t + \partial_{xxx}) u 
&= (\partial_t + \partial_{xxx}) (v_x + v^2) \\
&= \partial_x (\partial_t + \partial_{xxx}) v + 2v (\partial_t + \partial_{xxx}) v + 6 v_x v_{xx} \\
&= \partial_x (6 v^2 v_x) + 12 v^3 v_x + 6 v_x v_{xx}\\
&= \partial_x (6 v^2 v_x + 3 v^4 + 3 v_x^2) \\
&= \partial_x (3 u^2).
\end{align*}

The Miura transform acts roughly like a derivative, and in particular maps 
$H^s(\R)$ to $H^{s-1}(\R)$:

\begin{lemma}\label{continuity}  For any $0 \leq s < 1/2$ and $r > 0$, 
the Miura transform $M$ is Lipschitz continuous from the ball $\{ v_0 \in H^s(\R): \|v\|_{H^s(\R)} \leq r\}$ to $H^{s-1}(\R)$.   
\end{lemma}

\begin{proof}
The continuity is clear for the linear portion $v_x$ of the transform.  To obtain continuity for the quadratic portion $v^2$ we use the bilinear estimate
\be{hs}
\| v w \|_{H^{s-1}_x} \lesssim \| v w \|_{L^1_x} \lesssim \| v\|_{L^2_x} \| w \|_{L^2_x} \lesssim \| v \|_{H^{s}_x} \| w \|_{H^{s}_x}
\end{equation}
coming from H\"older and Sobolev embedding.
\end{proof}

From this we expect to use Theorem \ref{defocus-mkdv} to obtain Theorem \ref{kdv-illp}.  Unfortunately, due to low frequency issues, the Miura transform is not bilipschitz from $H^s$ to $H^{s-1}$, and so one must do a little computation (cf. the corresponding argument in \cite{kpv:counter}).

Let $0 \leq s < 1/4$, and let $0 < \delta \lesssim \eps \ll 1$ and $0 < T_0$ be given.  By the results in Section \ref{crude-sec} we can find smooth global solutions $\phi^{\langle a \rangle}$, $\phi^{\langle a' \rangle}$ to the mKdV equation \eqref{mkdv} such that
\bas
\| \phi^{\langle a \rangle}(0) \|_{H^s_x}+ \| \phi^{\langle a' \rangle}(0) \|_{H^s_x} &\lesssim \eps\\
\| \phi^{\langle a \rangle}(0) - \phi^{\langle a' \rangle}(0) \|_{H^s_x} &\lesssim \delta \\
\sup_{0 \leq t < T} \| \phi^{\langle a \rangle}(t) - \phi^{\langle a' \rangle}(t) \|_{H^s_x} &\gtrsim \eps
\end{align*}
for some small $0 < T \ll \log N/\lambda^3 \ll T_0$, 
where $N$, $\lambda$ are as in Section \ref{crude-sec}. 
In particular, $\lambda = N^{(1/4 -s)/(1/2+s)}$ tends to $\infty$
as $N\to\infty$.
Let $\psi^{\langle a \rangle} := M(\phi^{\langle a \rangle})$.  Then by Lemma \ref{continuity} we have
$$
\| \psi^{\langle a \rangle}(0) \|_{H^{s-1}_x}, \| \psi^{\langle a' \rangle}(0) \|_{H^{s-1}_x} \lesssim \eps
$$
and
$$
\| \psi^{\langle a \rangle}(0) - \psi^{\langle a' \rangle}(0) \|_{H^{s-1}_x} \lesssim \delta.$$
To finish the argument we would like to show
$$
\sup_{0 \leq t < T} \| \psi^{\langle a \rangle}(t) - \psi^{\langle a' \rangle}(t) \|_{H^{s-1}_x} \gtrsim \eps,$$
but this not quite automatic because $M$ is not bilipschitz at low frequencies.  To get around this we shall need more explicit control on the $\psi^{\langle a \rangle}$, using the details of the construction in Section \ref{crude-sec}.  

By \eqref{u-decohere} we may find $0 \leq t < T$ such that the NLS solutions $u^{\langle a\rangle}$, $u^{\langle a' \rangle}$ constructed in Section \ref{pseudo-sec} obey
\be{u-split}
\| u^{\langle a \rangle}(t\lambda^3) - u^{\langle a' \rangle}(t\lambda^3) \|_{H^5_x} \gtrsim \eps.
\end{equation}
Fix this $t$.  By the arguments in Section \ref{crude-sec}, this implies that
$$
\| \phi^{\langle a \rangle}(t) - \phi^{\langle a' \rangle}(t) \|_{H^s_x} \gtrsim \eps.$$
From \eqref{phi-bound} 
we have
$$
\| \phi^{\langle a \rangle}(t)\|_{H^s_x} + \|\phi^{\langle a' \rangle}(t) \|_{H^s_x} \lesssim \eps.$$
From \eqref{ee} and Lemma \ref{continuity} we have
$$ \| \psi^{\langle a \rangle}(t) - M(\phi^{[a]})(t) \|_{H^{s-1}_x} 
=  \| M(\phi^{\langle a\rangle})(t)  - M(\phi^{[a]})(t) \|_{H^{s-1}_x}
\lesssim \eps N^{-3/4-s} \ll \eps.$$
It will thus suffice to show that
$$ \| M(\phi^{[a]})(t) - M(\phi^{[a']})(t) \|_{H^{s-1}_x} \gtrsim \eps.$$
From \eqref{hs} we have
$$ \| \phi^{[a]}(t)^2 - \phi^{[a']}(t)^2 \|_{H^{s-1}_x} \lesssim
\| \phi^{[a]}(t) - \phi^{[a']}(t) \|_{H^s_x} 
\| \phi^{[a]}(t) + \phi^{[a']}(t) \|_{H^s_x}
\lesssim \eps^2 \ll \eps,$$
so it suffices to show 
$$ \| \partial_x (\phi^{[a]}(t) - \phi^{[a']}(t)) \|_{H^{s-1}_x} \gtrsim \eps.$$
We recall that for any $a$, $\phi^{[a]}(t)$ has the explicit form
$$ \phi^{[a]}(t,x) 
:= \lambda \sqrt{\frac{2}{3N}} \Re~  e^{iN\lambda x} e^{iN^3 \lambda^3 t} u^{\langle aw \rangle}(\frac{1}{\lambda^3 t+1},\frac{\lambda x}{\lambda^3 t+1}).$$
We subtract $\phi^{[a']}$ from $\phi^{[a]}$ and differentiate in $x$.  The worst term arises when the derivative hits $e^{iN\lambda x}$; by using 
\eqref{u-split}, Lemma \ref{modulate}, and \eqref{lam-def} we see that the $H^{s-1}_x$ norm of this term is $\gtrsim \eps$.  If the derivative hits the factor
$u^{\langle aw \rangle} - u^{\langle a'w \rangle}$, the resulting term
is much smaller; indeed, \eqref{uaw-bound}, Lemma \ref{modulate}, and \eqref{lam-def} imply that the $H^{s-1}$ norm of this term is 
$O(\eps \lambda/N)$, which equals $O(\eps N^{-(1/4+2s)/(1/2+s)})$.  
Since we are assuming that $s\ge 0$, this is $o(1)\cdot\eps$
as $N\to\infty$.
The claim then follows from  the triangle inequality.  This completes the proof of Theorem \ref{kdv-illp}.

\section{A generalized Miura transform}\label{general-sec}

We saw in the previous section how the Miura transform can convert ill-posedness for mKdV in $H^s$ to ill-posedness for KdV in $H^{s-1}$.  One might hope to also use this transform to convert the well-posedness for mKdV in $H^{1/4}$ (from Theorem \ref{mkdv-lwp}, or Section \ref{mkdv-lwp-sec}) to well-posedness for KdV at $H^{-3/4}$. However a difficulty arises because the Miura transform $u = v_x + v^2$ is not invertible\footnote{For instance, if $v\in {\mathcal S}$ is real-valued 
then $\int u$ must be non-negative.
  Furthermore, if $u$ lies in the range of $M$, then the Schr\"odinger operator $-\frac{d}{dx}^2 + u = (\frac{d}{dx}+v)(\frac{d}{dx}+v)^*$ cannot have any negative eigenvalues. See also the paper \cite{AKS:miura} of Ablowitz, Kruskal and Segur in which the range of the Miura transform is described.}.

On the other hand, for high frequencies $|\xi| \gg 1$ the derivative operator $v \to v_x$ is invertible, and the lower order term $v^2$ is negligible (as can be seen by the amount of surplus regularity in \eqref{hs}).  So it seems the Miura transform would be invertible if we could omit low frequency errors.  

Motivated by this, we define a generalized Miura transform $M: H^{1/4}_x \times H^1_x \to H^{-3/4}_x$ by
\begin{equation}
\label{GenMiura}
 M(v,w) := v_x + v^2 + w,
\end{equation}
where $v$ and $w$ are complex valued.  (The $H^1$ regularity of the error is not particularly special -- anything between $H^{1/4}$ and $H^{5/4}$ will suffice for the argument below).

\begin{lemma}\label{miura}  The transform $M: H^{1/4}_x \times H^1_x \to H^{-3/4}_x$ is locally Lipschitz.  Also, for any $A > 0$ there exists a Lipschitz transform 
$W_A: H^{-3/4}_x \to H^{1/4}_x \times H^{1/4}_x$ such that $M \circ W_A$ is the identity on the ball $B_A := \{ u \in H^{-3/4}_x: \|u\|_{H^{-3/4}} \leq A \}$.
\end{lemma}

\begin{proof}
The continuity of $M$ is immediate from Lemma \ref{continuity}.  Now to construct the inverse map $W_A$.  Fix $A$, and let $P$ be a smooth Fourier projection to the region $|\xi| \gtrsim C_A$ for some large $C_A$ depending on $A$.  

We need to construct $v \in H^{1/4}_x$ and $w \in H^1_x$ such that $v_x + v^2 + w = u$.  We begin by constructing $v$.  Observe that for $C_A$, $C'_A$ large enough, the map 
$$ v \mapsto \partial_x^{-1} P(u - v^2)$$
is a contraction on the ball $\{ v \in H^{1/4}_x: \|v\|_{H^{1/4}} \leq C'_A A \}$.  Indeed we have (cf. \eqref{hs})
\bas 
\| \partial_x^{-1} P(u - v^2) \|_{H^{1/4}_x} &\lesssim \| u \|_{H^{-3/4}_x} + C_A^{-3/4} \| v^2 \|_{L^2_x} \\
&\lesssim A + C_A^{-3/4} \| v \|_{L^4_x}^2 \\
&\lesssim A + C_A^{-1/8} \| v \|_{H^{1/4}_x}^2,
\end{align*}
so the above map maps the ball to itself, and the contraction property can also be obtained by a similar argument.  Thus we can construct a $v$ in this ball such that
$$ v_x = P(u - v^2).$$
If one then sets $w := u - v_x - v^2 = (1-P) (u-v^2)$ and $W_A(u) := (v,w)$, we see that we have constructed a map with the desired properties (the Lipschitz behavior following from similar estimates to the above).
\end{proof}
 
Of course, the modified Miura transform $M(v,w)$ no longer transforms
mKdV to KdV.  However, if $(v,w)$ are smooth solutions to the mKdV-like system
\begin{equation}\label{mkdv-system}
\begin{split}
v_t + v_{xxx} &= 6 (v^2 + w) v_x \\
w_t + w_{xxx} &= 6 (v^2 + w) w_x \\
v(x,0) &= v_0(x)\\
w(x,0) &= w_0(x)
\end{split}
\end{equation}
then the function $u(t) := M(v(t),w(t))$ will satisfy \eqref{kdv} with $u_0 := M(v_0,w_0)$.  Indeed, we have{\footnote{Can this algebraic identity be explained within the AKNS \cite{AKNS} framework?}}
\begin{align*}
(\partial_t + \partial_{xxx}) u 
&= (\partial_t + \partial_{xxx}) (v_x + v^2 + w) \\
&= \partial_x (\partial_t + \partial_{xxx}) v + 2v (\partial_t + \partial_{xxx}) v + 6 v_x v_{xx} + (\partial_t + \partial_{xxx}) w\\
&= \partial_x (6 (v^2 + w) v_x) + 12 v (v^2 + w) v_x + 6 v_x v_{xx}
+ 6 (v^2 + w) w_x \\
&= \partial_x (6 v^2 v_x + 3 v^4 + 3 v_x^2 + 3 w^2 + 6 wv_x + 6 w v^2) \\
&= \partial_x (3 u^2).
\end{align*}

The transform \eqref{GenMiura} is an extension of Gardner's generalization \cite{MiuraKdVSurvey} of the Miura transform:  
If $a$, $b$ are constants and $y$ satisfies 
$$ y_t + y_{xxx} = 6 (a^2 y^2 + b y) y_x$$
then $u := a y_x + a y^2 + by$ satisfies \eqref{kdv}.  (Indeed, one applies the generalized Miura transform with $v := ay$ and $w := by$). 

\begin{proposition}\label{iteration}
The Cauchy problem \eqref{mkdv-system} is locally well-posed in the space $H^{1/4}_x \times H^1_x$.
\end{proposition}

\begin{proof}
Intuitively, the system \eqref{mkdv-system} is a hybrid of the mKdV equation at $H^{1/4}_x$ and the KdV equation at $H^1_x$, and so the results should follow from the arguments in \cite{kpv:gkdv}.  As the arguments below show, this will indeed be the case.

We recall the $X$ norm defined in \eqref{Xnorm}.  For technical reasons involving the fractional Leibnitz rule it will be convenient to replace this norm with the augmented Besov-type norm 
\begin{equation}
\label{Xstar}
 \| v \|_{X^*} := \| v \|_X + (\sum_{j \in \Z} \| Q_j v \|_X^2)^{1/2}
\end{equation}
where for each integer $j$, $P_j$ is a standard Littlewood-Paley projection in space to frequencies $|\xi| \lesssim 2^j$, and $Q_j := P_j - P_{j-1}$.  Since the $Q_j$ are essentially orthogonal in $H^s_x$, it is easy to see that the key energy estimate \eqref{energy} continues to hold when $X$ is replaced with $X^*$.

We shall iterate $(v,w)$ in the norm
\begin{equation}
\label{Xstarstar}
\| (v,w)\|_{X^{**}} :=  \| v \|_{X^*} + \| D^{3/4} w \|_{X^*}
\end{equation}
on the slab $[0,T] \times \R$ for some sufficiently small $T$. 
From the augmented version of \eqref{energy} and \eqref{mkdv-system} (noting that $D^{3/4}$ commutes with the Airy operator $\partial_t + \partial_{xxx}$) we have
$$ \| v \|_{X^*} + \|D^{3/4} w\|_{X^*} \lesssim \|v(0)\|_{H^{1/4}_x} + \| w(0) \|_{H^1_x} + \| D^{1/4}((v^2+w) v_x) \|_{L^1_t(L^2_x)} + 
\| D((v^2+w) w_x) \|_{L^1_t(L^2_x)}.$$
By a H\"older in time we can estimate the $L^1_t(L^2_x)$ norm by the $L^2_t(L^2_x)$ norm, gaining a power of $T^{1/2}$.  If $T$ is sufficiently small, we may use a standard contraction mapping argument to obtain local well-posedness provided that we show the 
trilinear and bilinear estimates
\begin{align}
\| D^{1/4}(v v' v''_x) \|_{L^2_x(L^2_t)} &\lesssim \| v \|_{X^*} \| v'\|_{X^*} \| v''\|_{X^*} \label{m1}\\
\| D(v v' w_x) \|_{L^2_x(L^2_t)} &\lesssim \| v \|_{X^*} \| v'\|_{X^*} \| D^{3/4} w\|_{X^*} \label{m3}\\
\| D(w w'_x) \|_{L^2_x(L^2_t)} &\lesssim \| D^{3/4} w \|_{X^*} \| D^{3/4} w\|_{X^*}\label{m4}\\
\| D^{1/4}(w v_x) \|_{L^2_x(L^2_t)} &\lesssim \| D^{3/4} w \|_{X^*} \| v\|_{X^*} \label{m2}.
\end{align}

To prove these estimates we first argue informally.  Heuristically, the worst terms should arise when the powers of $D$ fall on the roughest function.  Ignoring all other terms, we are left with proving the estimates
$$
\| v v' D^{1/4} u_x \|_{L^2_x(L^2_t)} \lesssim \| v \|_{X^*} \|v' \|_{X^*} \| u \|_{X^*}$$
and
$$
\| w D^{1/4} u_x \|_{L^2_x(L^2_t)} \lesssim \| D^{3/4} w \|_{X^*} \| u \|_{X^*}.$$
But these follow from H\"older, after estimating $v$, $v'$ in $L^4_x(L^\infty_t)$, $w$ in $L^2_x(L^\infty_t)$, and $D^{1/4} u_x$ in $L^\infty_x(L^2_t)$.

To argue more rigorously it is easiest to use Littlewood-Paley  decomposition\footnote{One can also proceed using the fractional Leibnitz rule but one has to be careful because of all the $L^\infty$ norm.  See \cite{kpv:gkdv} for a further discussion of this issue.} $1 = \sum_j Q_j$.  For sake of illustration we prove \eqref{m3}; the other estimates are similar.

We can expand the left-hand side of \eqref{m3} as
$$ \| \sum_{j_1, j_2, j_3} D( (Q_{j_1} v) (Q_{j_2} v') (Q_{j_3} w_x) ) \|_{L^2_t(L^2_x)}.$$
Let us first consider the contribution of the terms when $j_1, j_2 < j_3 - 10$.  In this case, the summands have Fourier transform in the region $|\xi| \sim 2^{j_3}$ and so by orthogonality we can estimate the previous by
$$ (\sum_{j_3} \| 2^{j_3} \sum_{j_1, j_2 < j_3 - 10} (Q_{j_1} v) (Q_{j_2} v') (Q_{j_3} w_x) \|_{L^2_t(L^2_x)}^2)^{1/2}.$$
We can collapse the summations to rewrite this as
$$ (\sum_{j_3} \| 2^{j_3} (P_{j_3-10} v) (P_{j_3-10} v') (Q_{j_3} w_x) \|_{L^2_t(L^2_x)}^2)^{1/2}.$$
Using H\"older as indicated in the non-rigorous argument, and observing the $P_j$, $Q_j$ are all bounded on $X^*$, we can estimate this by
$$ (\sum_{j_3} (2^{j_3} \| v \|_{X^*} \| v' \|_{X^*} 2^{-j_3} \| D^{3/4} Q_{j_3} \|_X)^2)^{1/2}$$
which is acceptable by the definition of $X^*$.

The terms when $j_2, j_3 < j_1 - 10$ or when $j_1, j_3 < j_2 - 10$ can be treated similarly (in fact one gets some additional exponential gains due to the derivative $\partial_x$ falling on a relatively low frequency).  Together, these three cases cover all the ``high-low'' interactions when one of the $j$'s is much larger than the other two.  The remaining terms can be grouped into several ``high-high'' interactions in which two of the $j$'s are comparable to each other, and the third is comparable or smaller.  A typical such group of interactions is
$$
\| \sum_{j_3} \sum_{j_3-10 \leq j_2 \leq j_3+10} \sum_{j_1 \leq j_3 + 10}
D( (Q_{j_1} v) (Q_{j_2} v') (Q_{j_3} w_x) ) \|_{L^2_t(L^2_x)}.$$
Collapsing the $j_1$ summation and using the triangle inequality, we can bound this by
$$
\sum_{j_3} \sum_{j_3-10 \leq j_2 \leq j_3+10}
\| D( (P_{j_3+10} v) (Q_{j_2} v') (Q_{j_3} w_x) ) \|_{L^2_t(L^2_x)}.$$
The expression inside the norm has frequency $\lesssim 2^{j_3}$, so we can estimate the $D$ by a $2^{j_3}$.  Applying H\"older as before, we can then estimate this by
$$
\sum_{j_3} \sum_{j_3-10 \leq j_2 \leq j_3+10} 2^{j_3}
\| v \|_{X^*} \| Q_{j_2} v' \|_{X} 2^{-j_3} \| D^{3/4} Q_{j_3} \|_X,$$
and the claim follows from Cauchy-Schwarz and the definition of the $X^*$ norms.  The other terms can be treated similarly (and in some cases one even gets some additional exponential gains).  We omit the details. 
\end{proof}

We combine Lemma \ref{miura} and Proposition \ref{iteration} to prove
Theorem \ref{kdv-endpoint}.

Let $\phi_m \in \Schwarz $, the Schwarz class, for all $ 1 \leq m \in
\Z$. Assume $\phi_m \rightarrow \phi$ in $H^{-3/4}_x$. Then there
exists $A>0$ such that $\| \phi \|_{H^{-3/4}_x} \leq A $ and we may
assume that $\|\phi_m \|_{H^{-3/4}_x} \leq A $. By Lemma \ref{miura},
we have that 
$W_A ( \phi_m ) = (v_m (0) , w_m (0)) \in H^{1/4}_x \times H^1_x $, and
there exists $B>0$ such that $\| (v_m (0), w_m (0)) \|_{H^{1/4}_x
  \times H^2_x } \leq B$. The Lipschitz continuity of $W_A$ implies
$(v_m (0) , w_m (0) ) \rightarrow (v(0), w(0) ) = W_A ( \phi ).$

Proposition \ref{iteration} implies there exists $T = T(B) > 0$ such
that for all $1 \leq m \in \Z$, there exists a uniquely defined
continuous map from 
\begin{equation*}
H^{1/4}_x \times H^1_x \longmapsto X^{**} \subset C^0 ([0,T]; H^{1/4}_x \times
H^1_x ) 
\end{equation*}
taking $(v_m (0) , w_m (0))$ to $(v_m , w_m )$, a solution of
the initial value problem for the modified KdV system. Furthermore,
the smoothness properties of the data persist during the evolution:
$(v_m (t), w_m (t) ) \in H^k_x \times H^k_x $ if $(v_m (0 ) , w_m (0)
) \in H^k_x \times H^k_x $. For $t \in [0,T]$, define $u_m (t) = M
(v_m (t), w_m (t) )$. Then $u_m$ is such that $\phi_m \longmapsto u_m
(t)$ is a smooth well understood KdV evolution, by the explicit
calculation following \eqref{mkdv-system}.

\begin{remark}
Uniqueness is known to hold in the setting of smooth KdV
solutions. Therefore, the preceding procedure for the construction of KdV
solutions (by sending $\phi \rightarrow W_A (\phi) = (v, w)$ then
evolving the mKdV system to $(v(t), w(t))$ and applying the
generalized Miura transform $M(v(t), w(t)) = u(t)$) is independent of
the parameter $A$.
\end{remark}

The continuity properties of the map $W_A$, the mKdV system
data-to-solution map, the time evolution of the mKdV system and the
generalized Miura transform combine to imply: The sequence $u_m$
converges to some limit $u$ in the $C^0 ([0,T]; H^{-3/4} )$ norm, 
and, in this way defines a locally Lipschitz continuous data-to-solution
map $H^{-3/4} \ni \phi \longmapsto u \in C^0 ([0,T], H^{-3/4} )$ for
the KdV initial value problem.


The justification that $u$ is a weak solution of KdV requires that we
show that $\varphi u^2  \in L^1_{xt}$ where $\varphi$ is a test function.
Thus, it suffices to show that $u \in L^2_{xt, loc}$. Lemmas \ref{continuity}
and \ref{miura} imply that $u$ has essentially the same regularity 
properties as $v_x$, where $v \in X^*$. Since $v \in X^*$, we know that 
$v \in X$ and (see \eqref{Xnorm}) therefore $\partial_x D^{\frac{1}{4}} v 
\in L^\infty_x L^2_t .$ Thus, $D^{\frac{1}{4}} u \in L^\infty_x L^2_t$
and $u \in L^2_{xt,loc}$.

\begin{remark}
The construction above does not identify a Banach space of functions
of spacetime contained in $C^0 ([0,T]; H^{-3/4})$ in which we may assert
uniqueness of solutions of the KdV initial value problem. In
contrast, the solutions constructed in \cite{kpv:kdv} for $H^s_x$ data
with $s> -3/4$ are known to be unique in the Bourgain space $X_{s,
  1/2+}$. A  subset of $C^0 ([0,T]; H^{-3/4}_x)$ in which such a
uniqueness property does hold is the generalized Miura image
$M([{\mbox{solutions of system \eqref{mkdv-system}}}] \cap X^{**} )$. 
However, this is a rather unsatisfying uniqueness criterion because it is not
easy to test whether a 
function
$u$ lies in this class. A more
satisfactory resolution of these issues would be to prove uniqueness
of KdV solutions evolving from $H^s_x$ data in the space $C^0 ([0,T];
H^s_x)$. This remains open for $-3/4 \leq s < 0$ but has been
established \cite{ZhouUnique} by Zhou for $s \geq 0$.

\end{remark}

\section{Proofs for the periodic case} \label{periodic-sec}

To prove Theorem~\ref{periodic-thm}, we begin with the simplest
equation, the defocusing nonlinear Schr\"odinger equation.
A family of explicit solutions is
\[
u_{N,a}(t,x) = ae^{i(Nx+N^2t-|a|^2t)},
\]
where $a\in\complex$ and $N$ is any positive integer.
Suppose that $s<0$.
Since $\|u_{N,a}(0,\cdot)\|_{H^s}\sim N^s$,
we choose $a=a(N)=N^{-s}\alpha$.
Comparing two solutions 
$u_{N,a}$ and $u_{N,a'}$
with $a = ,\alpha N^{-s}$ and $a' = \alpha'  N^{-s}$,
the $H^s$ norm of their difference at $t=0$ is
$O(|\alpha-\alpha'|)$, 
while there exist some $0<t\le T(N,\alpha,\alpha')
= C|\alpha-\alpha'|^{-2} N^{2s}$ for which the $H^s$ norm
of their difference is $\gtrsim|\alpha|+|\alpha'|$.
If $s<0$ then for any fixed $\alpha,\alpha'$, $T(N,\alpha,\alpha')$
may be made arbitrarily small, by choosing $N$ sufficiently large.
This implies illposedness in $H^s$.
Exactly this argument was given by Burq, G\'erard and Tzvetkov
\cite{BGT:NLSsphere}

\medskip 
Consider next the real mKdV equation.
Fix an exponent $s\in(-1,1/2)$.
We construct solutions $u=u_{N,\beta}$
of the form
\begin{align}
u(t,x) & = \sum'_k b_k e^{ik\psi}
\label{sumsolution}
\\
\psi(t,x) & = Nx+N^3t+\sigma t
\\
b_1&=\beta N^{-s}
\end{align}
where $\sigma,b_k$ are real numbers,
$b_{-k}\equiv b_k$ for all $k$,
$\beta>0$,
and the notation $\displaystyle\sum'$ means that the
sum extends over all {\em odd} integers.

Formally, such a function satisfies mKdV if and only if
the numerical coefficients satisfy the system of equations
\begin{equation} \label{odesystem}
\left\{
\begin{aligned}
\sigma &= b_1^{-1} 6N\sum'_{k_1+k_2+k_3=1} b_{k_1}b_{k_2}b_{k_3}
\\
b_k &= \frac{6N}{N^3(1-k^2)+\sigma}
\sum'_{k_1+k_2+k_3=k} b_{k_1}b_{k_2}b_{k_3}\ \text{ for all $|k|>1$.}
\end{aligned}\right.
\end{equation}
Here the sums extend over all three-tuples of odd integers satisfying
the stated relations.

We claim that, for any $s\in(-1,1/2)$,
for any $\delta>0$ and any $\beta\in[\delta,1]$,
for all sufficiently large $N$ there exist $\sigma,\{b_k\}$
satisfying
\begin{align*}
\sigma&\sim N^{1-2s}
\\
|b_k|&\le AN^{-s}N^{-\eta\cdot(|k|-1)}
\ \text{for all } |k|>1,
\end{align*}
where $\eta,A\in\R^+$ depend only on $\delta,s$.
In particular, for large $N$, the dominant terms in the sum
defining $u$ are those with $k=\pm 1$.
The restriction $s>-1$ means that $N^{1-2s}\ll N^3$,
hence the term $\sigma$ appearing in the denominator of the equation for $b_k$
for $|k|\ne 1$ is negligible.
These estimates imply that the series \eqref{sumsolution}
converges rapidly for large $N$, and that $\|u(0,x)\|_{H^s}\sim \beta$.

Because $1-2s>0$,
The factor $e^{i\sigma t}$ multiplying
$e^{i(Nx+N^3t)}$ causes solutions to become out of phase
within a timespan $\lesssim |\beta-\beta'|^{-2} N^{2s-1}$.
Because $1-2s>0$, this tends to zero as $N\to\infty$
so long as $\beta,\beta'$ are fixed.
Thus the equation is illposed for this range of exponents.

Existence of, and the bounds for, solutions of \eqref{odesystem}
for large $N$ can be proved by a straightforward application of
the contraction mapping principle; details are left to the reader.
This completes the discussion of periodic real mKdV.

\medskip
The illposedness for real KdV is obtained from that for real mKdV
via the Miura transform, as in Section~\ref{miura-sec}.

\begin{remark}
It is not clear whether one can create global periodic solutions to mKdV of the type in Sections \ref{energy-sec}, \ref{approx-sec}, mainly because one does not have the crucial decay in $t$.  (In the periodic case such decay is impossible, as can be seen e.g. by $L^2$ norm conservation).
\end{remark}


\end{document}